\documentclass{amsart}

\usepackage{latexsym}
\usepackage{amsmath}
\usepackage{amssymb}

\newtheorem{theorem}{Theorem}[subsection]
\newtheorem{lemma}[theorem]{Lemma}
\newtheorem{corollary}[theorem]{Corollary}

\theoremstyle{definition}
\newtheorem{definition}[theorem]{Definition}
\newtheorem{example}[theorem]{Example}

\theoremstyle{remark}
\newtheorem{remark}[theorem]{Remark}
\newtheorem{question}[theorem]{Question}

\numberwithin{equation}{section}

\newcommand{\G}{\mathcal{G}}

\DeclareMathOperator{\tor}{Tor}

\newcommand{\comp}[1]{#1^c}
\newcommand{\numbc}{\mbox{\#$\operatorname{comp}$}}
\newcommand{\reg}{\operatorname{reg}}

\newcommand{\pdim}{\operatorname{pd}}
\newcommand{\lcm}{\operatorname{lcm}}
\newcommand{\compGS}{\comp{G_{S}}}
\newcommand{\N}{\mathbb{N}}
\newcommand{\Z}{\mathbb{Z}}

\newcommand{\lk}{\operatorname{conn}}
\newcommand{\rlk}{\overline{\operatorname{conn}}}

\newcommand{\I}{\mathcal{I}}
\newcommand{\F}{\mathcal{F}}
\newcommand{\A}{\mathcal{A}}
\newcommand{\K}{\mathcal{K}}

\begin{document}

\title[Resolutions of monomial ideals via facet ideals]
{Resolutions of square-free monomial ideals via facet ideals: a 
survey}

\author{Huy T\`ai H\`a}
\address{Department of Mathematics\\ Tulane University 
\\ 6823 St. Charles Ave. \\ New Orleans, LA 70118, USA}
\email{tha@tulane.edu}
\urladdr{http://www.math.tulane.edu/$\sim$tai/}

\author{Adam Van Tuyl}
\address{Department of Mathematical Sciences \\
Lakehead University \\
Thunder Bay, ON P7B 5E1, Canada}
\email{avantuyl@sleet.lakeheadu.ca}
\urladdr{http://flash.lakeheadu.ca/$\sim$avantuyl/}

\date{Version: \today}

\begin{abstract} 
We survey some recent results on the minimal graded free resolution of a
square-free monomial ideal.  The theme uniting these results is the 
point-of-view
that the generators of a monomial ideal correspond to the maximal faces 
(the facets)
of a simplicial complex $\Delta$.  
\end{abstract}
\maketitle

%%%%%%%%%%%%%%%%%%%%%%%%%%%%%%%%%%%%%%%%%%%%%%%%%%%%%%%
\setlength{\parskip}{3pt}
\section{Introduction}
Let $I$ be a {\it monomial ideal} in a polynomial ring $R = 
k[x_1,\ldots,x_n]$.  
Associated to $I$ is a {\it minimal
graded free resolution} of the form
\[
0 \rightarrow \bigoplus_j R(-j)^{\beta_{l,j}(I)}
\rightarrow \bigoplus_j R(-j)^{\beta_{l-1,j}(I)}
\rightarrow \cdots
\rightarrow  \bigoplus_j R(-j)^{\beta_{0,j}(I)}
\rightarrow I \rightarrow 0
\]
where $l \leq n$ and $R(-j)$ is the free $R$-module obtained by 
shifting
the degrees of $R$ by $j$.  The number $\beta_{i,j}(I)$, the $ij$th
{\it graded Betti number} of $I$, equals the number of minimal 
generators
of degree $j$ in the $i$th syzygy module of $I$.

A classical topic in commutative algebra is to understand how 
the graded Betti numbers in the minimal free resolution 
of a monomial ideal depend upon the generators of the ideal.
This problem  continues to inspire current research;  we refer the 
reader to 
Miller and Sturmfels' book \cite{MS} and Villarreal's book
\cite{V2} for background, descriptions of various approaches 
to the problem, and many relevant references to ongoing work.

This paper surveys a new perspective on the problem of understanding
the resolution of a monomial ideal that has appeared since \cite{MS}.
The new point-of-view relates the
graded Betti numbers of monomial ideals to combinatorial objects 
described by the 
generators of the monomial ideal.  More precisely, we know from work
of Fr\"oberg \cite{Fr1} that the study of graded Betti numbers of 
monomial ideals can be reduced
to understanding the case that $I$ is generated by {\it square-free 
monomials}. 
We then have a bijection between the sets
\[
\left\{
\mbox{simplicial complexes $\Delta$ on $n$ vertices}
\right\}
\leftrightarrow
\left\{
\begin{array}{c}
\mbox{square-free monomial} \\
\mbox{ideals $I \subseteq R = k[x_1,\ldots,x_n]$}
\end{array}
\right\} 
\]
given by
\[\Delta \stackrel{1-1}{\longleftrightarrow} \I(\Delta) = \left.\left
\langle \prod_{x \in F} x ~\right|~ 
F ~\mbox{is a {\it facet} (a maximal face) of 
$\Delta$}\right\rangle.\]
The ideal $\I(\Delta)$ is called the {\it facet ideal} of $\Delta$\footnote[1]{Since first submitting this paper, some authors
\cite{HHTZ,HVT2,VTV} have taken the point of view that
the monomial generators of a square-free monomial ideal correspond to the edges of a hypergraph
$\mathcal{H}$ (sometimes called a clutter).  A {\it hypergraph}
is pair $\mathcal{H} = (\mathcal{X},\mathcal{E})$ where $\mathcal{X}
= \{x_1,\ldots,x_n\}$ are the vertices, and $\mathcal{E} = \{E_1,\ldots,
E_s\}$ is edges with $E_i \subseteq \mathcal{X}$.  The
edge ideal of $\mathcal{H}$ is
$\I(\mathcal{H}) = (\{ \prod_{x \in E} x ~|~ \mbox{$E$ an edge of 
$\mathcal{H}$}\}).$  One then wishes to describe the resolution of 
$\I(\mathcal{H})$
in terms of the properties of $\mathcal{H}$.  This approach
is similar to the one described above, the only difference being
whether the combinatorial object is viewed as a
 simplicial complex or a hypergraph. }.
We 
then wish
to describe the minimal graded free resolution of $\I(\Delta)$ in terms 
of the combinatorial
data of $\Delta$.  The novelty of this approach is
to view the generators of a square-free monomial ideal as the maximal 
faces of
the simplicial complex.  This contrasts with the usual Stanley-Reisner 
correspondence
which associates to a simplicial complex $\Delta$ the square-free 
monomial ideal 
$I_{\Delta}$ generated by
the minimal {\it non-faces} of $\Delta$. 

When all the facets of $\Delta$ have dimension one, $\Delta$ can be 
viewed as a simple 
graph $G$ (a graph with no loops or multiple edges) on $n$ vertices.  
We shall usually write 
$\I(G)$ in this case
for $\I(\Delta)$, and we call $\I(G)$ the {\it edge ideal} of $G$.
Observe that $\I(G)$ is generated by square-free quadratic monomials
and is the first non-trivial case of a square-free monomial
ideal.
Historically, edge ideals were introduced by Villarreal in \cite{V1} 
before facet ideals.
Facet ideals, which can be seen as a generalization of edge ideals, 
were introduced later by 
Faridi in \cite{faridi:2002}.  To learn more about the properties of 
edge ideals,
one should see \cite{B1,B2,DFT,FH,FV,HHZ2, HHZ,HHZ1,R,S,SVV, SS,Su,V,V2}; further 
properties of facet ideals can be found
in the sequels \cite{faridi:2004,faridi:2005,faridi:2005b} to Faridi's 
paper 
cited above.

Eliahou and Villarreal \cite{EV} provided one of the first 
examples showing that the numbers $\beta_{i,j}(\I(\Delta))$ could be 
described in terms of the 
combinatorial data of $\Delta$.  Specifically, it was shown
that the number $\beta_{1,3}(\I(G))$ could be computed
from the degrees of the vertices of $G$ and the number of triangles of 
$G$.  
Although Fr\"oberg's paper \cite{Fr} predates the notion of an edge 
ideal,  Fr\"oberg
demonstrated that a facet ideal point-of-view could provide information 
about
the resolution of a monomial ideal; using  the language of edge ideals, 
the
ideal $\I(G)$ has a linear resolution if and only if $G^c$, the 
complement
of $G$, is a {\it chordal} graph.  Since 2003, there has been a flurry of 
results
on the resolutions of facet (and edge) ideals.  We mention
in particular Corso and Nagel \cite{CN},
Eisenbud, Green, Hulek, Popescu \cite{eghp}, Horwitz \cite{Hor},
Jacques \cite{J}, Jacques and Katzman \cite{JK}, Katzman \cite{K}, 
Visscher \cite{Viss},
Zheng \cite{Z}, the second author and Roth \cite{RVT}, and the two 
authors \cite{HVT,HVT2}.

The goal of this survey paper is two-fold.  Our first goal is to 
provide a
summary of the state-of-the-art on the resolutions of facet and edge 
ideals.  
This is a relatively new area, and we hope to gather together most of 
the relevant results 
currently available in the literature. We will not concentrate on 
proving these results, 
but rather, we try to develop a unified and systematic perspective to 
these studies. From 
time to time, we shall give the reader a flavor of the proofs by 
sketching out 
important ideas. 

As we see it, there are two main themes in the literature.  The first 
theme 
is to study the structure of 
the resolution of edge and facet ideals. For example, Katzman \cite{K} 
showed that the minimal 
free resolution of edge ideals of certain graphs depends upon char$(k)$,
the characteristic of $k$;
Eisenbud, et. al. \cite{eghp} characterized property $N_{2,p}$ of edge 
ideals via the cycle 
structure 
of complement graphs; and Zheng \cite{Z} calculated the regularity of 
the edge ideal of a forest from the 
number of disconnected edges in the graph.  More recently,
Corso and Nagel \cite{CN}, Horwitz \cite{Hor} and Visscher \cite{Viss}
have been interested in describing the maps in the resolutions
of edge ideals. 
The second theme is to give explicit 
computations for 
various graded Betti numbers. For instance, Jacques \cite{J} provided 
formulas for 
the graded Betti numbers of edge ideals of special classes of graphs 
including
cycles and complete graphs; Roth and the second author \cite{RVT} 
calculated the graded 
Betti numbers in the linear strand of edge ideals of graphs having no 
minimal cycle of 
length 4;
Zheng \cite{Z} computed the graded Betti numbers in the linear strand 
of facet ideals of 
simplicial forests; and the two authors \cite{HVT,HVT2} derived 
recursive-type formulas for the
graded Betti numbers of the edge ideal of a graph in terms of its 
subgraphs.

To present the known results, we have divided the results into two
main categories based upon the techniques used in their proves.  In 
Section 3,
we describe what results have been obtained using the theory of reduced 
simplicial
homology.  In particular, we describe results that one can obtain 
using Hochster's formula
and its variant, the formula of Eagon and Reiner.  In Section 4 we 
describe the results
on resolutions of facet ideals which
rely on the notion of
a splittable ideal (first introduced by Eliahou and Kervaire 
\cite{EK}).

The second goal of this paper is based upon the belief that a good 
survey
should also inspire future research on a topic.  To this end,  we 
provide
in Section 5
a collection of open questions that we would like answered.  The 
questions have been
grouped into four broad categories.  The first category is related
to building a dictionary between the graded Betti numbers 
$\beta_{i,j}(\I(\Delta))$
and the combinatorial data of $\Delta$.  The second category discusses 
questions on
how the characteristic of the ground field affects the numbers 
$\beta_{i,j}(\I(\Delta))$.
The third group of questions deal with other homological invariants such 
as
regularity and projective dimension.  The final category is a series
of questions based upon the authors' paper \cite{HVT} and the notion of 
splittable
ideals.    

It is our hope that the reader will be convinced that  there is an 
advantage of viewing 
square-free monomial
ideals as the facet ideal of a simplicial complex.  However, we do wish 
to point out  some
natural limitations of this approach.  G. Reisner \cite{Re} showed that 
the
square-free monomial ideal 
\begin{eqnarray*}
\I(\Delta) &=& 
(x_1x_2x_3,x_1x_2x_4,x_1x_3x_5,x_1x_4x_6,x_1x_5x_6,x_2x_3x_6,\\
&& x_2x_4x_5,x_2x_5x_6,x_3x_4x_5,
x_3x_4x_6)
\end{eqnarray*}
in $R =k[x_1,\ldots,x_6]$ is an ideal whose resolution depends upon 
char$(k)$.  Specifically,
$\I(\Delta)$ has a linear resolution except in the case that char$(k) = 
2$.
This example illustrates that knowing the combinatorial object 
$\Delta$ will not always be enough to 
understand the minimal graded free resolution (see also \cite{K} for an 
example 
involving an edge ideal).  
So although we cannot expect a theory that relates
the data of every simplicial complex $\Delta$ to the resolution 
of $\I(\Delta)$,  we do not feel that this diminishes the importance
of this new perspective.  As shown in this paper, 
the combinatorial data is enough to describe the entire (or part of 
the) resolution in
many interesting cases, and we are of the belief that many more 
interesting cases
remain to be discovered using this approach. 
%%%%%%%%%%%%%%%%%%%%%%%%%%%%%%%%%%%%%%%%%%%%%%%%%%%%%%%

\section{Preliminaries}

To make this paper as self-contained as possible, we have included the 
relevant
definitions about graphs, simplicial complexes, and minimal 
resolutions.

\subsection{Graphs}
Let $G$ denote a graph with vertex set $V_G$ and edge set $E_G$.
We shall say that $G$ is {\it simple}
if $G$ has no loops or multiple
edges. A simple graph need not be connected.
We shall let $\numbc(G)$ denote
the number of connected components of $G$. 
The {\it degree} of a vertex
$x \in V_G$, denoted by $\deg x$, is the number of edges incident to 
$x$.
 
If $S \subseteq V_G$, then the {\it induced subgraph} of $G$ on the
vertex set $S$, denoted by $G_{S}$, is the
subgraph of $G$ whose vertex set is $S$ and whose edge set consists of 
edges of $G$ connecting two vertices in $S$.  The {\it complement} of a 
graph $G$, denoted
by $\comp{G}$, is the graph whose vertex set is the
same as $G$, but whose edge set is defined by
the rule: $\{x_i,x_j\} \in E_{G^c}$ if and only
if $\{x_i,x_j\} \not\in E_G$.

We write $xy$ as a shorthand for the edge $\{x,y\}$. A {\it cycle} of 
length $q$ in a graph $G$ is 
a sequence of edges $\{ e_1 = x_1x_2, e_2 = x_2x_3, \dots, e_{q-1} = 
x_{q-1}x_q, e_q = x_qx_1 \}$ 
in $G$ (where $x_i \not= x_j$ for $i \not= j$). 
We use $(x_1x_2 \cdots x_qx_1)$ to denote a cycle of length $q$ with 
vertices $x_1, \dots, x_q$. 
We sometimes also use $C_q$ to refer to a cycle of length $q$. We say 
that a 
cycle $C = (x_1x_2 \cdots x_qx_1)$ of $G$ has a {\it chord} if there 
exists some 
$j \not\equiv i+1 (\text{mod} \ q)$ such that $x_ix_j$ is an edge of 
$G$. We call a 
cycle $C$ in $G$ a {\it minimal cycle} if $C$ has length at least 4 and 
contains no chord. 

A {\it forest} is any graph with no cycles; a {\it tree} is a connected 
forest.
The {\it wheel} $W_n$ is the graph
obtained by adding a vertex $z$ to $C_n$ and then
adjoining an edge between $z$ and every vertex in $V_{C_n}$.
Note that $W_n$ has $n+1$ vertices.
The {\it complete graph} on $n$ vertices, denoted $\K_n$,
is the graph with
the property that for all $x_i, x_j \in V_{\K_n}$ with $i \neq j$,
the edge $x_ix_j\in E_{\K_n}$.  The {\it complete bipartite graph},
denoted $\K_{n,m}$, is the graph with vertex set $V_G = 
\{x_1,\ldots,x_n,y_1,
\ldots,y_m\}$ and edge set $E_G = \{x_iy_j ~|~ 1 \leq i \leq n,~~
1 \leq j \leq m\}$.  We shall write $c_n(G),w_n(G),k_n(G),$ and 
$k_{n,m}(G)$
for the number of induced subgraphs of $G$ isomorphic to
$C_n,W_n,\K_n$, and $\K_{n,m}$, respectively.

\subsection{Simplicial complexes}  
A {\it simplicial complex} $\Delta$ on a vertex set
$V_\Delta$ is a collection of subsets
of $V_\Delta$ such that for all $x \in V_\Delta$, $\{x\} \in \Delta$,
and for each $F \in \Delta$, if $G \subseteq F$, then
$G \in \Delta$.  Note that $\emptyset \in \Delta$, except in
the case that $\Delta =\{ \}$ is the void complex (see \cite[Definition 
1.4]{MS}).

An element of a simplicial complex $\Delta$ is called a {\it face} of 
$\Delta$.
The {\it dimension} of a face $F$ of $\Delta$, denoted $\dim F$, is 
defined to
 be $|F|-1$, where $|F|$ denotes the number of vertices in $F$. The 
{\it dimension}
 of $\Delta$, denoted by $\dim \Delta$, is defined to be the maximal 
dimension of a
 face in $\Delta$. The maximal faces of $\Delta$ under inclusion are 
called {\it facets}
 of $\Delta$. If all facets of $\Delta$ have the same dimension $d$, 
then $\Delta$ is
said to be {\it pure $d$-dimensional}.
 
We usually denote the simplicial complex $\Delta$ with facets $F_1, 
\dots, F_q$
by
$$\Delta = \langle F_1, \dots, F_q \rangle;$$
here, the set $\F(\Delta) = \{ F_1, \dots, F_q \}$ is often referred to 
as the
{\it facet set} of $\Delta$. There is a one-to-one correspondence 
between simplicial 
complexes and their facet sets. If we realize a simplicial complex 
$\Delta$ by its 
facet set $\F(\Delta)$, then a graph $G$ can be thought of as a 
simplicial complex 
with the facet set being its edge set. Suppose $F$ is a facet of 
$\Delta$, say $F = F_q$, then
 we denote by $\Delta \backslash F$ the simplicial complex 
obtained by removing
 $F$ from $\F(\Delta)$, i.e., $\Delta \backslash F = \langle F_1, \dots, 
F_{q-1} \rangle.$
Throughout the paper, by a {\it subcomplex} of a simplicial complex 
$\Delta$, we  
shall mean a simplicial complex whose facet set is a subset of the 
facet set of
$\Delta$. If $\Delta'$ is a subcomplex of $\Delta$, then we denote by
$\Delta \backslash \Delta'$ the simplicial complex obtained from 
$\Delta$ by
 removing from its facet set all facets of $\Delta'$.
 
We say that two facets $F$ and $G$ of $\Delta$ are {\it connected} if 
there exists a
{\it chain} of facets of $\Delta$, $F = F_0, F_1, \dots, F_m = G$, such 
that
$F_i \cap F_{i+1} \not= \emptyset$ for any $i = 0, \dots, m-1$. The 
simplicial
complex $\Delta$ is said to be {\it connected} if any two of its facets 
are connected.

A facet $F$ of a simplicial complex $\Delta$ is a {\it leaf}
of $\Delta$ if either $F$ is the only facet of $\Delta$, or there 
exists a facet  
$G$ in $\Delta$, $G \not= F$, such that $F \cap F' \subseteq F \cap G$ 
for every  
facet $F' \in \Delta, F' \not= F$. 
It follows from \cite[Remark 2.3]{faridi:2002}
that if $F$ is a leaf of $\Delta$, then $F$ must contain a vertex that 
does not
belong to any other facet of the complex (the converse is not true 
though). 
A simplicial complex $\Delta$ is
called a {\it tree} if $\Delta$ is connected and every nonempty 
connected
subcomplex of $\Delta$ (including $\Delta$ itself) has a leaf. We call 
$\Delta$ 
a {\it forest} if every connected component of $\Delta$ is a tree.

If $\Delta$ is a simplicial complex over a vertex set $V_\Delta = 
\{x_1, \dots, x_n\}$, 
then we can 
associate to $\Delta$ two ideals in the polynomial ring $R = k[x_1, 
\dots, x_n]$, the facet ideal and the 
Stanley-Reisner ideal (by abuse of notation, we shall use $x_i$s to 
denote both the vertices 
of $\Delta$ and the variables in the polynomial ring). For a face $F$ 
of $\Delta$, we write 
$x^F$ to denote the monomial $\prod_{x \in F} x$ in $R$. The {\it facet 
ideal} of $\Delta$ 
is defined to be 
$$\I(\Delta) = \langle x^F ~|~ F \in \F(\Delta) \rangle \subseteq R,$$ 
and the {\it Stanley-Reisner ideal} of $\Delta$ is defined to be 
$$I_\Delta = \langle x^F ~|~ F \subseteq V_{\Delta}, F \not\in \Delta 
\rangle
 \subseteq R.$$

Finally, we can associate to any graph a simplicial complex. 
If $G$ is any graph, then the {\it clique complex} of $G$ is the 
simplicial
complex $\Delta(G)$ where $F = \{x_{i_1},\ldots,x_{i_j}\} \in
\Delta(G)$ if and only if $G_{F}$ is a complete graph.
Note that if $G$ is a simple graph with edge ideal $\I(G)$, then
the ideal $\I(G)$ is generated by square-free monomials.
So $\I(G)$ is also the Stanley-Reisner ideal of a simplicial complex
$\Delta$, that is, $\I(G) = I_{\Delta}.$ Specifically,
$\I(G)$ is the Stanley-Reisner ideal associated to
the clique complex $\Delta = \Delta(\comp{G})$
constructed from the complement of $G$.

\subsection{Minimal resolutions}

Let $R = k[x_1,\ldots,x_n]$.
If $\cdots \rightarrow \mathcal{F}_1 \rightarrow \mathcal{F}_0 \rightarrow I \rightarrow 
0$
is a minimal graded free resolution of $I$, then 
$\cdots \rightarrow \mathcal{F}_1 \rightarrow \mathcal{F}_0 \rightarrow R \rightarrow R/I 
\rightarrow 0$
is a minimal graded free resolution of $R/I$.  From this observation
it follows that $\beta_{i,j}(R/I) = \beta_{i-1,j}(I)$ for all $i,j \geq 
0$. Here, we will adopt the 
convention that $\beta_{-1,0}(I) = 1$ and $\beta_{-1,j}(I) = 0$ for any 
$j > 0$.
 
The {\it regularity} of $I$,  denoted $\reg(I)$, is defined to be
\[\reg(I) := \max\{j-i ~|~ \beta_{i,j}(I) \neq 0\}.\]
The {\it projective dimension} of $I$, denoted $\pdim(I)$, is defined 
to be
\[\pdim(I):= \max\{i ~|~ \beta_{i,j}(I) \neq 0 \}.\]
These invariants 
measure the ``size'' of the  minimal graded free resolution.
 
If $d$ is the smallest degree of a generator of an ideal $I$, then the 
Betti numbers
$\beta_{i,i+d}(I)$ form the
so-called {\it linear strand} of $I$ (see \cite{eghp,hi}).  An  ideal 
$I$ generated
by elements all of degree $d$  is said to have a {\it linear 
resolution}
if $\beta_{i,j}(I) = 0$ for all $j \neq i+d$, that is, the only nonzero 
graded
Betti numbers are those in the linear strand.
 
\begin{example} (The numbers $\beta_{0,j}(\I(\Delta))$)  
The number $\beta_{0,j}(\I(\Delta))$ is simply the number of 
generators of $\I(\Delta)$
of degree $j$.  From our construction of $\I(\Delta)$ it follows
that $\beta_{0,j}(\I(\Delta)) = $ number of facets of $\Delta$ of 
dimension $j-1$.  When
$\Delta = G$ is a simple graph, then  $\beta_{0,j}(\I(G))$ equals 
$|E_G|$, the number of edges
if $j=2$, and equals $0$ if $j \neq 2$.  Since the relation between 
the
$0$th graded Betti numbers and the combinatorics of $\Delta$ is well 
understood, 
we can restrict to studying the numbers $\beta_{i,j}(\I(\Delta))$ with 
$i \geq 1$.
\end{example}

\begin{example}\label{beta0j}
 (Disconnected simplicial complexes)  Suppose
that $\Delta$ is a simplicial complex that is the disjoint union
of two components, i.e. 
$\Delta = \Delta_1 \cup \Delta_2$ with $\Delta_1 \cap \Delta_2 = 
\emptyset$.
If $V_{\Delta_1} = \{x_1,\ldots,x_n\}$ and $V_{\Delta_2} = 
\{y_1,\ldots,y_m\}$,
then $\I(\Delta) = \I(\Delta_1) + \I(\Delta_2)$ and 
\[\frac{k[x_1,\ldots,x_n,y_1,\ldots,y_m]}{\I(\Delta)} \cong 
\frac{k[x_1,\ldots,x_n]}{\I(\Delta_1)} \otimes_k 
\frac{k[y_1,\ldots,y_m]}{\I(\Delta_2)}.\]
It then follows that the graded Betti numbers of $\I(\Delta)$ can be 
calculated
by finding the resolutions of $\I(\Delta_1)$ and $\I(\Delta_2)$ and 
then tensoring the two resolutions together (see \cite[Lemma 2.1]{JK} 
for details).
The upshot is that to study the numbers $\beta_{i,j}(\I(\Delta))$ one
can make the harmless assumption that $\Delta$ is connected.
\end{example}

%%%%%%%%%%%%%%%%%%%%%%%%%%%%%%%%%%%%%%%%%%%%%%%%%%%%%%%

\section{Results via reduced simplicial homology}

Given a square-free monomial ideal $I$, one can associate to $I$ two 
simplicial complexes.
The first is the simplicial complex $\Delta$ with $I = \I(\Delta)$;
the second is the simplicial complex $\Delta'$ with $I = I_{\Delta'}$
via the standard Stanley-Reisner correspondence.  The starting point 
for
most of the known results in the literature about the graded Betti 
numbers of a facet ideal 
$\I(\Delta)$ is to first describe the associated simplicial complex 
$\Delta'$ 
with the property that $\I(\Delta) = I_{\Delta'}$.  Then one appeals 
to
results such as  the formulas of Hochster and Eagon-Reiner
to describe the graded Betti numbers of $I_{\Delta'}$ in terms of the
reduced simplicial homology of $\Delta'$.  Finally, one translates
these results back in terms of the simplicial complex $\Delta$.  
This approach has proven extremely useful in describing the 
graded Betti numbers of edge ideals of a simple graphs $G$ due, in 
part, to the fact that
we know that $\Delta' = \Delta(G^c)$, the clique complex of the 
corresponding complement graph.  
In this section, we shall describe what formulas and results 
about facet ideals have been derived via this approach.

\subsection{Hochster's formula and Eagon-Reiner's formula}
Hochster's and Eagon-Reiner's formulas describe the graded Betti 
numbers of
a square-free monomial ideal $I$ in terms of the simplicial complex 
$\Delta$
where $I = I_{\Delta}$.  We recall these methods below.

Let $\Delta$ be a simplicial complex over a vertex set $V_\Delta = 
\{x_1, \dots, x_n\}$, and let $I_\Delta$ be
the Stanley-Reisner ideal of $\Delta$. We can view $R$ as an 
$\N^n$-graded $k$-algebra and $I_\Delta$ as a $\Z^n$-graded $R$-module. For a 
monomial $m$ of $R = k[x_1,\ldots,x_n]$ 
having multidegree $\alpha
\in \N^n$, we define
\[\tor_i^R(I_{\Delta},k)_m := \tor_i^R(I_{\Delta},k)_{\alpha}.\]
Also, for $W \subseteq V_{\Delta}$ we write $\Delta_{W} = \{ F \in 
\Delta
~|~ F \subseteq W\}$ for the restriction of $\Delta$ to $W$. 
It follows
that $\Delta_W$ is a simplicial complex on $W$.
If $m$ is a square-free
monomial ideal, and $W$ denotes the set of variables that divide $m$, then we will 
write
$|m|$ for $\Delta_{W}$.
Hochster \cite{Ho} provided the following formula to 
calculate $\beta_{i,j}(I_\Delta)$.

\begin{theorem}[Hochster's Formula]\label{prop: hochster} 
Let $\Delta$ be a simplicial complex on the vertex set
$V_\Delta = \{x_1,\ldots,x_n\}$ and let $m$ be a monomial of $R$.
If $m$ is square-free, then
\[\dim_k \tor_i^R(I_{\Delta},k)_m =
\dim_k \widetilde{H}_{\deg(m)-i-2}(|m|,k)\]
where $\widetilde{H}_j(|m|,k)$ denotes the $j$th reduced homology
of $|m|$.  If $m$ is not square-free,
then $\tor_i^R(I_{\Delta},k)_m$ vanishes.
In particular
\[\beta_{i,j}(I_{\Delta}) = 
\sum_{\footnotesize
\begin{array}{c} 
\deg m =j, ~~m ~\mbox{is square-free}
\end{array}}
\dim_k \widetilde{H}_{j-i-2}(|m|,k)~~\mbox{for all $i,j \geq 0$.}\]
\end{theorem}

Eagon and Reiner \cite{ER} introduced a variant of Hochster's formula
that uses Alexander duality.  The {\it Alexander dual} of a simplicial 
complex $\Delta$ is given by
\[\Delta^{\vee} := \{\{x_1,\ldots,x_n\} \backslash F ~|~ F \not\in 
\Delta\}.\]
Moreover, if $F \in \Delta$, then the {\it link} of $F$ in $\Delta$, 
denoted by Link$_{\Delta} F$, is the
simplicial complex defined by
\[\operatorname{Link}_{\Delta}F := \{G \in \Delta ~|~ G \cup F \in 
\Delta
~\mbox{and}~ G \cap F = \emptyset\}.\]

\begin{theorem}[Eagon-Reiner's Formula] \label{Eagon-Reiner}
Let $\Delta$ be a simplicial complex. Then
\[\beta_{i,j}(I_{\Delta}) = \sum_{F \in \Delta^{\vee}, ~|F|=n-j}
\dim_k \tilde{H}_{i-1}(\operatorname{Link}_{\Delta^{\vee}} F, k) 
~~\mbox{for all $i,j \geq 0$.}\]
\end{theorem}

As shown by Jacques \cite{J}, Jacques
and Katzman \cite{JK}, and Katzman \cite{K}, it sometimes is easier 
to compute the graded Betti numbers of edge ideals via the 
Eagon-Reiner
reformulation of Hochster's formula because it is easier to deal
with the reduced homology of the simplicial complexes 
Link$_{\Delta^{\vee}} F$.

\begin{remark} (Characteristic dependence of Betti numbers)
From the conclusions of Theorems \ref{prop: hochster} and 
\ref{Eagon-Reiner} 
it becomes clear that the characteristic of the field simply cannot be 
ignored since
the characteristic may introduce some nonzero torsion.  An example
due to Reisner can be found in the introduction.  
Katzman \cite{K} shows
that even if $\I(G)$ is an edge ideal, the graded Betti numbers are 
not independent of char$(k)$.  
See Section 5 for some questions about this topic.
\end{remark}

\subsection{Graded Betti numbers of edge ideals}
Recall that when $G$ is a simple graph we have $\I(G) = I_{\Delta}$, 
where $\Delta = \Delta(G^c)$ is
the clique complex of the complementary graph $G^c$.
This observation allows us to use Hochster's and Eagon-Reiner's 
formulas to compute the graded 
Betti numbers of $\I(G)$ (a similar formula is found in \cite{K}).

\begin{theorem}\label{prop: edgeidealbetti}
Let $G$ be a simple graph with edge ideal $\I(G)$.  Then
\[\beta_{i,j}(\I(G)) = \sum_{
\footnotesize
\begin{array}{c} 
S \subseteq V_G, ~~
|S| = j
\end{array}}
 \dim_k \widetilde{H}_{j-i-2}(\Delta(\comp{G_{S}}),k)
~~\mbox{for all $i,j\geq 0$.}\]
\end{theorem}

\begin{proof} We sketch out the main idea of the proof.
Let $\Delta = \Delta(G^c)$ be the simplicial complex defined by 
$\I(G)$.
It follows from Proposition \ref{prop: hochster} that
\[\beta_{i,j}(\I(G)) = \beta_{i,j}(I_{\Delta})
= \sum_{\footnotesize
\begin{array}{c} 
m \in M_j, ~~m ~\mbox{is square-free}
\end{array}}
\dim_k \widetilde{H}_{j-i-2}(|m|,k)\]
where $M_j$ consists of all the monomials of degree $j$ in $R$.
Since $\deg m = j$ and $m$ is square-free,
the variables that divide $m$ give a subset
$S \subseteq V_G$ of size $j$.  
Let $G_{S}$ denote the induced subgraph of $G$ on this vertex
set $S$, and let $\comp{G_{S}}$ denote its complement.  To finish the 
proof, 
it is enough to note that $|m|$, the restriction 
of $\Delta(G^c)$ to $S$, and 
$\Delta(\comp{G_{S}})$ are the same simplicial complex.
\end{proof}

\begin{example} \label{beta25}
We illustrate how to apply the above theorem to compute 
$\beta_{2,5}(\I(G))$.
Below are all the graphs $H$ on $5$ vertices with 
$\dim_k \widetilde{H}_{1}(\Delta(H^c),k) > 0$:
\begin{picture}(400,45)(0,-5)
\put(0,15){$H_1 =$}
\put(30,0){\circle*{2}}
\put(50,0){\circle*{2}}
\put(30,0){\line(1,0){20}}
\put(30,20){\line(1,0){20}}
\put(30,20){\circle*{2}}
\put(30,20){\line(1,1){10}}
\put(50,20){\circle*{2}}
\put(40,30){\line(1,-1){10}}
\put(40,30){\circle*{2}}

\put(60,15){$H_2 =$}
\put(90,0){\circle*{2}}
\put(110,0){\circle*{2}}
\put(90,0){\line(1,0){20}}
\put(90,20){\circle*{2}}
\put(90,20){\line(1,1){10}}
\put(110,20){\circle*{2}}
\put(100,30){\line(1,-1){10}}
\put(100,30){\circle*{2}}

\put(120,15){$H_3 =$}
\put(150,0){\circle*{2}}
\put(150,0){\line(0,1){20}}
\put(170,0){\circle*{2}}
\put(170,0){\line(0,1){20}}
\put(150,0){\line(1,0){20}}
\put(150,20){\circle*{2}}
\put(150,20){\line(1,1){10}}
\put(170,20){\circle*{2}}
\put(160,30){\line(1,-1){10}}
\put(160,30){\circle*{2}}

\put(180,15){$H_4 =$}
\put(210,0){\circle*{2}}
\put(210,0){\line(0,1){20}}
\put(230,0){\circle*{2}}
\put(230,0){\line(0,1){20}}
\put(210,20){\circle*{2}}
\put(210,20){\line(1,1){10}}
\put(230,20){\circle*{2}}
\put(220,30){\line(1,-1){10}}
\put(220,30){\circle*{2}}

\put(240,15){$H_5 =$}
\put(270,0){\circle*{2}}
\put(270,0){\line(1,0){20}}
\put(290,0){\circle*{2}}
\put(290,0){\line(0,1){20}}
\put(270,20){\circle*{2}}
\put(270,20){\line(1,0){20}}
\put(270,20){\line(1,1){10}}
\put(290,20){\circle*{2}}
\put(280,30){\line(1,-1){10}}
\put(280,30){\circle*{2}}

\put(300,15){$H_6 =$}
\put(330,0){\circle*{2}}
\put(330,0){\line(1,0){20}}
\put(330,0){\line(1,1){20}}
\put(350,0){\circle*{2}}
\put(350,0){\line(0,1){20}}
\put(330,20){\circle*{2}}
\put(330,20){\line(1,0){20}}
\put(330,20){\line(1,1){10}}
\put(350,20){\circle*{2}}
\put(340,30){\line(1,-1){10}}
\put(340,30){\circle*{2}}
\end{picture}

For $i = 2,\ldots,6$,  $\dim_k \widetilde{H}_1(\Delta(H^c_i),k) = 1$, 
and
$\dim_k \widetilde{H}_1(\Delta(H^c_1),k) = 2$ (these numbers are 
independent
of $k$).  If we let $h_i(\Gamma)$ denote the number
of induced subgraphs of the graph $\Gamma$ isomorphic to $H_i$, 
then by the Theorem \ref{prop: edgeidealbetti}
we have
\[\beta_{2,5}(\I(G)) = 2h_1(G) + h_2(G) + h_3(G) + h_4(G) 
+h_5(G) + h_6(G).\] 
\end{example}

As evident by the previous example, the formula of Theorem \ref{prop: 
edgeidealbetti} is
difficult to apply since one has to compute 
the dimensions of all the homology groups 
$\widetilde{H}_{j-i-2}(\Delta(\comp{G_{S}}),k)$ as
$S$ varies over all subsets of $V_G$ of size $j$.
It is also not immediately clear how this formula relates to 
combinatorial data, like
degrees of the vertices or the number of cliques, associated to
$G$.  However, in special cases, we can still tease out interesting 
conclusions (as we will
see below).

We now show how Theorem \ref{prop: edgeidealbetti} (and other tools) 
can be
used to give exact formulas for some of the graded Betti numbers of 
$\I(G)$
in terms of $G$.
We begin by observing that we can restrict our search for 
graded Betti numbers to particular ranges.

\begin{theorem}\label{thm: betti edge range}
Let $G$ be a simple graph with edge ideal $\I(G)$.
If $\beta_{i,j}(\I(G)) \neq 0$, then  $i+2 \leq j \leq 2(i+1)$.
\end{theorem}

\begin{proof}
Since $\I(G)$ is generated by quadrics, $\beta_{i,j}(\I(G)) = 0$ if $j 
< i+2$.  
By using the Taylor resolution it can also be seen that
$\beta_{i,j}(\I(G)) = 0$ for all $j > 2(i+1)$ (see Katzman \cite{K} for
details).  For more on the Taylor resolution, and a generalization
of its construction,
see Herzog \cite{Herz}.
\end{proof}

For the extremal values of $j$, that is, $j = i+2$ or $j=2(i+1)$,
we can compute $\beta_{i,j}(\I(G))$ in terms of data from $G$ for each 
$i$.
Observe that these numbers are independent of char$(k)$ since they only 
depend upon
the graph $G$.

\begin{theorem}\label{prop: linearstrand}
Let $G$ be a simple graph with edge ideal $\I(G)$.
Then for all $i \geq 0$
\begin{eqnarray*}
\beta_{i,i+2}(\I(G)) &= &
\sum_{\footnotesize
\begin{array}{c} 
S \subseteq V_G, ~~|S| = i+2
\end{array}}
 \left(\numbc(\comp{G_{S}})-1\right) ~\mbox{and}\\
\beta_{i,2(i+1)}(\I(G)) &= &\left|\left\{ H ~\left|~ 
\begin{array}{l}
\mbox{$H$ is a induced subgraph of $G$} \\
\mbox{consisting of $i+1$ disjoint edges}
\end{array}\right\}\right.\right|.
\end{eqnarray*}
\end{theorem}

\begin{proof}
The formula for $\beta_{i,i+2}(\I(G))$ is given in \cite[Proposition 
2.1]{RVT}; it is a
consequence of evaluating the formula of Theorem \ref{prop: 
edgeidealbetti} at $j=i+2$ and
using the fact that $\dim_k \tilde{H}_0(\Gamma,k)+ 1$ is the number of 
connected components
of $\Gamma$.   The formula
for $\beta_{i,2(i+1)}(\I(G))$ comes from \cite[Lemma 2.2]{K} and relies 
on the Taylor resolution.
\end{proof}

\begin{remark} From Theorem \ref{prop: linearstrand}
the length of the linear strand is given by
\[\ell = \max\{i ~|~ \mbox{there exists $S \subseteq V_G$ with 
$|S|=i+2$ and $\numbc(\comp{G_{S}}) > 1$}
\}.\]
\end{remark}

Theorem \ref{prop: linearstrand} gives us a means 
to compute the graded Betti numbers in the linear strand of $\I(G)$.
However, one is required to sum over all subgraphs of a certain size
which limits the usefulness of the result.

Roth and the second author \cite{RVT} showed that in
many cases one can find equivalent (and more easy to calculate) 
formulas
for the graded Betti numbers in the linear strand.
Their results are based upon the following decomposition  
for the formula for $\beta_{i,i+2}(\I(G))$:  
\begin{equation}
\label{splitequation}
\footnotesize
\beta_{i,i+2}(\I(G)) =
\sum_{\footnotesize
\begin{array}{c}
S \subseteq V_G,\\
|S| = i+2, \\
G^c_{S}\, \mbox{contains an}\\
\mbox{isolated vertex}\\
\end{array}}
 (\numbc(\compGS)-1)
+
\sum_{\footnotesize
\begin{array}{c}
S \subseteq V_G, \\
|S| = i+2,\\
G^c_{S}\, \mbox{contains no}\\
\mbox{isolated vertices}\\
\end{array}}
 (\numbc(\compGS)-1).
\end{equation}
Recall that  $k_i(G)$ is the number of induced subgraphs of $G$
isomorphic to $\K_i$.

\begin{theorem}\cite[Proposition 2.4]{RVT} \label{prop: 
linearstrand_no_c4}
Let $G$ be a simple graph with edge ideal $\I(G)$.  If $G$
has no minimal $4$-cycles, then
\[\beta_{i,i+2}(\I(G)) = \sum_{v \in V_G} \binom{\deg v}{i+1}
- k_{i+2}(G) 
~~\mbox{for all $i \geq 0$.}\]
Furthermore, the above formula holds for all simple graphs $G$ if $i=0$ 
or $1$.
\end{theorem}

\begin{proof}
The proof has two steps.  The first step is to show that the 
second
sum in (\ref{splitequation}) is $0$.  So, suppose that $S \subseteq 
V_G$
is such that $G_S^c$ has no isolated vertex.
We claim that $G_S^c$ cannot have two connected components.  By 
checking
all possible graphs on $2$ or $3$ vertices, this is easy to
see if $|S|=2,3$.
If $|S| \geq 4$ and if $G_S^c$ has at least two connected
components, then 
there must be at least two edges with each edge in a different 
connected component.  Let $S'$ be the
set of four vertices incident to these two edges.  Then $G_{S'}$ is 
a minimal $4$-cycle, contradicting our hypothesis.  Thus,
if $G_S^c$ has no isolated vertex, then $G_S^c$ is connected,
and thus makes no contribution to $\beta_{i,|S|}(\I(G))$.

The second step is to evaluate the first sum.  We first observe
that if $G_S^c$ has an isolated vertex, it can have at most
one connected component consisting of one or more edges.  If $G_S^c$
had two or more connected components having an edge, then by argument
similar to the one given above, this would imply that $G$ has a 
minimal
4-cycle.
Therefore, to count $\numbc(\compGS)-1$ we can simply count the number 
of
isolated vertices in $\compGS$, the
``$-1$'' term being taken care of by the component which is not a 
vertex;
we thus over count by one whenever
$\compGS$ consists completely of isolated vertices.
 
For any vertex $v$, the number of subsets $S$ of size $i+2$ containing
$v$ such that $v$ is an isolated vertex in $\compGS$ is
$\binom{\deg v}{i+1}$.
To take care of the over count, note that $\compGS$ consists of 
isolated
vertices
exactly when $G_{S}$ is a complete graph on $i+2$ vertices.  
Subtracting the
number of times this happens gives the formula above.
 
To compute $\beta_{i,i+2}(\I(G))$
when $i = 0$ or $1$, we need to count the number of
connected components of $\comp{G_{S}}$ when $|S| =2$ or $3$.
But for any simple graph on two or three vertices, at most
one connected component can be larger than a vertex.  The proof
is now the same as the one given above.\end{proof}

Since a forest has no cycles (and hence, no induced subgraphs isomorphic to 
$\K_j$ with
$j \geq 3$) we obtain: 

\begin{corollary}\cite[Corollary 2.6]{RVT} \label{foreststrand}
Let $G$ be a forest with edge ideal $\I(G)$.  Then
$\beta_{0,2}(\I(G)) = |E_G|$, and 
\[\beta_{i,i+2}(\I(G)) = \sum_{v \in V_G} \binom{\deg v}{i+1}
~~\mbox{for all $i \geq 1$.}\]
\end{corollary}

\begin{example}\label{beta1j}
(The numbers $\beta_{1,j}(\I(G))$)
The above results allow us to completely describe $\beta_{1,j}(\I(G))$ 
for
all $j$ and all graphs $G$:
\[
\beta_{1,j}(\I(G)) =
\left\{
\begin{array}{ll}
{\displaystyle \sum_{v \in V_G}\binom{\deg v}{2}- k_3(G)} & \mbox{if $j 
= 3$} \\
c_4(G^c) & \mbox{if $j = 4$} \\
0 & \mbox{otherwise.}
\end{array}\right.
\]
The fact that $\beta_{1,j}(\I(G)) = 0$ if $j \neq 3,4$ comes from 
Theorem \ref{thm: betti edge range}.
The formula for $\beta_{1,3}(\I(G))$ is just Theorem \ref{prop: 
linearstrand_no_c4}.
By Theorem \ref{prop: linearstrand} 
\[\beta_{1,4}(\I(G)) =|\{H ~|~ \mbox{$H$ is an induced subgraph of $G$ 
with exactly 2 disjoint edges}\}|.\]
But any induced subgraph $H$ that is exactly 2 disjoint edges 
corresponds to an induced cycle $C_4$ in
$G^c$.  So $\beta_{1,4}(\I(G))$ is the number of 4-cycles in $G^c$.
The formula for $\beta_{1,3}(\I(G))$ was first proved in 
\cite[Proposition 2.1]{EV}.
\end{example}

By a  careful analysis of the second sum
in (\ref{splitequation}), Roth and the second author derived
formulas for $\beta_{2,4}(\I(G))$ and $\beta_{3,5}(\I(G))$
for any simple graph $G$.  The formula for $\beta_{2,4}(\I(G))$ 
verifies a
conjecture found in \cite[Conjecture 2.4]{EV}.

\begin{theorem}\cite[Proposition 2.8]{RVT} \label{exact form}
Let $G$ be a simple graph with edge ideal $\I(G)$.  Then
\begin{eqnarray*}
\beta_{2,4}(\I(G)) &=& \sum_{v\in V_G} \binom{\deg v }{3} - k_4(G) + 
k_{2,2}(G)~~\mbox{and} \\
\beta_{3,5}(\I(G)) & = &\sum_{v \in V_G}\binom{\deg v}{4} - k_5(G) + 
k_{2,3}(G)
+ w_4(G) + d(G)
\end{eqnarray*}
where  $d(G)$ 
is the number of induced 
subgraphs of $G$ isomorphic to the graph:

\begin{picture}(-100,40)(-155,0)
\put(0,0){\circle*{5}}
\put(0,0){\line(1,0){30}}
\put(0,0){\line(0,1){30}}
\put(0,30){\circle*{5}}
\put(0,30){\line(1,0){30}}
\put(30,0){\circle*{5}}
\put(30,0){\line(0,1){30}}
\put(30,30){\circle*{5}}
\put(15,15){\circle*{5}}
\put(0,0){\line(1,1){30}}
\put(0,30){\line(1,-1){15}}
\end{picture}
\end{theorem}

\begin{example}\label{beta2j}
(The numbers $\beta_{2,j}(\I(G))$)  By Theorem \ref{thm: betti edge 
range}
the number $\beta_{2,j}(\I(G))$ 
is nonzero only if $j =4,5,$ or $6$.  The formula for 
$\beta_{2,4}(\I(G))$ 
can be found in the above theorem.  We computed $\beta_{2,5}(\I(G))$ in 
Example \ref{beta25}.  Finally,
one can compute $\beta_{2,6}(\I(G))$ by using Theorem \ref{prop: 
linearstrand}.  Notice that
the numbers $\beta_{2,j}(\I(G))$ are independent of char$(k)$ for all $j$.
\end{example}

When $G$ has a minimal $4$-cycle, we can compute upper and lower 
bounds
on the graded Betti numbers in the linear strand.

\begin{theorem} \cite[Proposition 3.1]{RVT} \label{prop: 
linearstrand_lower_bound}
Let $G$ be a simple graph with edge ideal $\I(G)$.
Suppose that $G$ has an induced $4$-cycle.  Then for all $i \geq 2$
{\small
\[\beta_{i,i+2}(\I(G)) \geq \sum_{v \in V_G} \binom{\deg v}{i+1}
- k_{i+2}(G) + k_{2,i}(G) + k_{3,i-1}(G) + 
\cdots + k_{\lfloor \frac{i+2}{2}\rfloor,\lceil 
\frac{i+2}{2}\rceil}(G).\]} 
\end{theorem}

\begin{proof}
The first sum in (\ref{splitequation}) is bounded below by  $ \sum_{v 
\in V_G}\binom{\deg v}{i+1}
- k_{i+2}(G)$, while 
the second sum is bounded below
by  $k_{2,i}(G) + k_{3,i-1}(G) + 
\cdots + k_{\lfloor \frac{i+2}{2}\rfloor,\lceil 
\frac{i+2}{2}\rceil}(G)$.  
\end{proof}

Recall that the elements of $\N^n$ can be given a total ordering
using the lexicographical order defined by $(a_1,\ldots,a_n)
> (b_1,\ldots,b_n)$ if $a_1 = b_1, \ldots, a_{i-1} = b_{i-1}$
but $a_i > b_i$.  This induces an ordering on the monomials
of $R$:  $x_1^{a_1}\cdots x_{n}^{a_n} >_{lex} x_1^{b_1}\cdots 
x_n^{b_n}$
if $(a_1,\ldots,a_n) > (b_1,\ldots,b_n)$.  A monomial
ideal $I$ is a {\it lex ideal} if for each $d \in \N$,
a basis for $I_d$ is the $\dim_k I_d$ largest monomials of degree $d$
with respect to the lexicographical ordering. 

\begin{theorem}\cite[Proposition 3.2]{RVT} \label{prop: upperbound}
Let $G$ be a simple graph with edge ideal $\I(G)$. If
$\{m_1,\ldots,m_{|E_G|}\}$ are the $|E_G|$ largest monomials
of degree 2 in $R$ with respect to the lexicographical 
ordering,  then
\[\beta_{i,i+2}(\I(G)) \leq \sum_{t=1}^{|E_G|} \binom {u_t -1}{i}\]
where $u_t$ is the largest index of a variable dividing $m_t$. 
\end{theorem}

\begin{proof}
One uses Eliahou-Kervaire's
formula \cite{EK} for the growth of the graded Betti numbers of stable 
(and hence
lex) ideals to bound $\beta_{i,i+2}(\I(G))$.
\end{proof}

We now turn our attention to the global behavior of the resolutions. 
The algebraic invariants and properties in which we shall be interested 
include the regularity,  the linear strand, and property $N_{2,p}$.

The first such result places a lower bound on the regularity
of any edge ideal.  Moreover, this bound is exact when
$G$ is a chordal graph.  
Zheng \cite{Z} introduced a notion of two edges being disconnected; precisely,
two edges $u_1v_1$ and $u_2v_2$ of a simple graph $G$ are {\it 
disconnected} \footnote[2]{Two edges that are disconnected
according to Zheng's definition are called $3$-disjoint by the two
authors in \cite{HVT2}.} if 
(a) the two edges do not share a common vertex, and 
(b) $u_1u_2,u_1v_2,v_1u_2,v_1v_2$ are not 
edges of $G$.  
Note that a pair of disconnected edges can belong to the
same same connected component of $G$.  Alternatively, if $d(x,y)$ 
denotes
the distance between the vertices $x$ and $y$, that is, the least 
length
of a path from $x$ to $y$, then Zheng's definition is equivalent to 
saying
that two edges $u_1v_1$ and $u_2v_2$ are disconnected if
 $d(u_1,u_2), d(u_1,v_2),d(v_1,u_2)$ and $d(v_1,v_2)$ are all at least 
$2$.

\begin{theorem}\cite[Theorem 6.5, Corollary 3.9]{HVT2}\label{tree reg}
Let $G$ be a graph with edge ideal $\I(G)$.  If $j$ is the maximal 
number of pairwise disconnected edges in $G$, then
\[\reg(\I(G)) \geq j+1.\]
If $G$ is a chordal graph, then the above inequality is
an equality.
\end{theorem}

\begin{remark}
The above result was first proved for the
case that $G$ is forest in \cite[Theorem 2.18]{Z}.
\end{remark}

Zheng \cite[Remark 2.19]{Z} points out that if $G = 
C_5$,
the $5$-cycle, then $\I(G)$
is an example where
$\reg(\I(G)) \neq $ maximal number of pairwise disconnected edges $+1$.  
However, for edge ideals, there is also an upper bound for
the regularity, using the matching number.

\begin{definition}
A {\bf matching} of $G$ is a set of pairwise disjoint edges.
The {\bf matching number} of $G$, denoted $\alpha'(G)$, is the largest size
of a maximal matching in $G$.
\end{definition}

\begin{theorem} \label{upperbound}
Let $G$ be a finite simple graph.  Then
$$\reg(R/\I(G))\leq \alpha'(G)$$
where $\alpha'(G)$ is the matching number of $G$.
\end{theorem}

\begin{proof} The proof (which can be found in \cite{HVT2})
is based upon the Taylor resolution of $\I(G)$.
\end{proof}

Recall that a
cycle $C$ is a {\it minimal cycle} if $C$ has length at least 4 and 
contains no chord. An 
ideal $I$ is said to satisfy {\it property $N_{2,p}$} for some $p \ge 
1$ if $I$ is 
generated by quadratics and its minimal free resolution is linear up to 
the $p$th step, 
i.e., $\beta_{i,j}(I) = 0$ for all $0 \le i \le p-1$ and $j > i+2$. 
Eisenbud, Green, Hulek and 
Popescu \cite{eghp} gave an interesting characterization for property 
$N_{2,p}$ for edge ideals
in terms of the minimal cycles of $G^c$.
We restate this result as follows:

\begin{theorem}\cite[Theorem 2.1]{eghp} \label{np property}
Let $G$ be a simple graph with edge ideal $\I(G)$. Then $\I(G)$ 
satisfies property $N_{2,p}$ with $p > 1$ 
if and only if every minimal cycle in $G^c$ has length $\ge p+3$.
\end{theorem}

\begin{proof} Hochster's formula (Theorem \ref{prop:  hochster}) 
is employed liberally throughout the proof of \cite{eghp}.
A careful study of the
reduced simplicial homology groups $\tilde{H}_i(|m|,k)$, where $|m|$ 
denotes the restriction
of $\Delta = \Delta(G^c)$ to the vertices corresponding to the 
variable dividing $m$, is
required in the proof.  
We will sketch out an alternative combinatorial proof in the next 
section (see Corollary \ref{N2P2}).
\end{proof}

Since a chordal graph is a graph that has no minimal cycles, the 
following result of Fr\"oberg \cite{Fr}
becomes a corollary of Theorem \ref{np property}.

\begin{corollary}\cite[Theorem 1]{Fr} \label{linear resolution} Let $G$ 
be a graph with edge ideal $\I(G)$.
Then $\I(G)$ has a linear resolution if and only if $G^c$ is a
chordal graph.
\end{corollary} 

\begin{remark}
Reisner's example given in the introduction shows that we cannot expect a 
purely combinatorial description of the simplicial complexes
$\Delta$ with the property that $\I(\Delta)$ has a linear resolution.
(For Reisner's example, the resolution of $\I(\Delta)$ is linear
if and only if the characteristic is not two).  The papers
of Bruns and Hibi \cite{BH,BH2} look at the question of when
one can combinatorially identify  simplicial complexes that must have a 
linear (or pure) resolution.
\end{remark}

\begin{example}\label{k1dbettinumbers}(The resolution of $\I(\K_{a,b})$)
Let $G = \K_{a,b}$ be a complete bipartite graph.  We write
the vertex set of $G$ as $V_G = \{x_1,\ldots,x_a,y_1,\ldots,y_b\}$
so that $E_G = \{x_iy_j ~|~ 1 \leq i \leq a, ~1\leq j \leq b\}$.
For all $a,b \geq 1$, the complement of $G$ is the
disjoint union of $\K_a$ and $\K_b$.  Since $G^c$
has no induced cycles of length $\geq 4$, the resolution
of $\I(G)$ is linear by Corollary \ref{linear resolution}.
 
Because $\comp{G} = \K_a \cup \K_b$,
for any $S \subseteq V_G$ with $|S| = i+2$, we have
\[
\numbc(\comp{G_{S}}) =
\left\{
\begin{array}{ll}
2 & \text{if}~S \cap\{x_1,\ldots,x_a\} \neq \emptyset
~\text{and}~S \cap\{y_1,\ldots,y_b\} \neq \emptyset \\
1 & \mbox{otherwise.}
\end{array}
\right.
\]
By Theorem \ref{prop: linearstrand}, to determine
$\beta_{i,i+2}(\I(G))$ it therefore suffices to
count the number of subsets $S \subseteq V_G$ with
$|S| = i+2$ and $\numbc(\comp{G_{S}}) =2$.
 
There are $\binom{a+b}{i+2}$ subsets of $V_G$ that contain
$i+2$ distinct vertices.  Furthermore, $\binom{a}{i+2}$
of these subsets must contain only vertices among
$\{x_1,\ldots,x_a\}$; similarly, $\binom{b}{i+2}$
of these subsets contain only vertices among $\{y_1,\ldots,y_b\}$.
It thus follows that
\[\beta_{i,i+2}(\I(\K_{a,b})) = \binom{a+b}{i+2} -
\binom{a}{i+2} - \binom{b}{i+2} ~~\text{for all $i \geq 0$}\]
since the expression on the right hand side counts the number
of subsets $S \subseteq V_G$ with $|S| = i+2$ and $S$ contains
at least one $x_i$ vertex and one $y_j$ vertex.

By adapting this proof, one can find formulas for the graded Betti
numbers of the edge ideals for the multipartite graph 
$\K_{d_1,\ldots,d_n}$
(see also \cite[Theorem 5.3.8]{J}).  In the recent paper
of Visscher \cite{Viss}, the the maps in the minimal free resolution of 
$\I(\K_{a,b})$ are also described.  Similar results can
also be found in the paper of Corso and Nagel \cite{CN} in
which they study the edge ideals of Ferrers graphs, a class
of graphs which includes all complete bipartite graphs.
\end{example}

%%%%%%%%%%%%%%%%%%%%%%%%%%%%%%%%%%%%%%%%%%%%%%%%%%%%%%

\section{Splittable monomial ideals}

In this section we present a new tool for investigating the graded 
Betti
numbers of facet (and edge) ideals.  This approach, which
uses the notion of a splittable ideal introduced by Eliahou and 
Kervaire
(see \cite{EK}), was first
explored by the authors \cite{HVT,HVT2}.
This notion is surprisingly strong in that a splittable monomial
approach produces results from which previously known results can
be deduced as corollaries.  As well, this tool gives us a means to 
generalize
results known about the resolutions of edge ideals to the more
general situation of facet ideals.

\subsection{Splittable ideals}
For a monomial ideal $I$, let $\G(I)$ denote the set of minimal 
monomial generators of $I$; 
this set is uniquely 
determined (cf. \cite[Lemma 1.2]{MS}). 

\begin{definition}\label{defn: split}
A monomial ideal $I$ is {\it splittable} if $I$ is the sum
of two nonzero monomial ideals $J$ and $K$, that is, $I = J+K$, such 
that
\begin{enumerate}
\item $\G(I)$ is the disjoint union of $\G(J)$ and $\G(K)$.
\item there is a {\it splitting function}
\begin{eqnarray*}
\G(J\cap K) &\rightarrow &\G(J) \times \G(K) \\
w & \mapsto & (\phi(w),\psi(w))
\end{eqnarray*}
satisfying
\begin{enumerate}
\item for all $w \in \G(J \cap K), ~~ w = \lcm(\phi(w),\psi(w))$.
\item for every subset $S \subset \G(J \cap K)$, both 
$\lcm(\phi(S))$ and $\lcm(\psi(S))$
strictly divide $\lcm(S)$.
\end{enumerate}
\end{enumerate}
If $J$ and $K$ satisfy the above properties, then
$I = J + K$ is a {\it splitting} of $I$.
\end{definition}

When $I = J + K$ is a splitting of a monomial
ideal $I$, then there is a relation between  $\beta_{i,j}(I)$ 
and the graded Betti numbers of the ``smaller'' ideals. 

\begin{theorem}[Eliahou-Kervaire \cite{EK}, Fatabbi \cite{Fa}]  
\label{prop: ekf}
Suppose that $I$ is a
splittable monomial ideal with splitting $I = J+K$.  Then
for all $i, j \geq 0$,
\[\beta_{i,j}(I) = \beta_{i,j}(J) + \beta_{i,j}(K) + 
\beta_{i-1,j}(J\cap K).\]
\end{theorem}

This theorem suggests an approach to the study of the numbers 
$\beta_{i,j}(\I(\Delta))$.  
Precisely, one wishes to find splittings of the ideal $\I(\Delta) = J + 
K$ such
that the ideals $J,K$, and $J \cap K$ are related to facet ideals of 
subcomplexes of $\Delta$.  
Theorem \ref{prop: ekf} then provides a recursive-type formula for the 
numbers
$\beta_{i,j}(\I(\Delta))$.  This is the general strategy of \cite{HVT} 
and  the results
of this approach are described below.  Note that
the formulas will not be recursive in general since
we may not be able to split the new facet ideals arising from $J,K$, 
and $J \cap K$.

\subsection{Splitting edges}

Let $G$ be a simple graph with edge ideal $\I(G)$ and  $e = uv \in 
E_G$.
If we set
\[ J = (uv) ~\mbox{and}~  K= \I(G\backslash e),\]
then $\I(G) = J + K$.  In general this 
may not be a splitting of $\I(G)$ because the second
condition of Definition \ref{defn: split} may not hold.
Thus, an edge $e = uv$ with the property that $\I(G) = (uv) + 
\I(G\backslash e)$
is a splitting is a special type of edge. We give such an edge the 
following name:

\begin{definition}An edge $e = uv$ is a 
{\it splitting edge} of $G$ if $J = (uv)$ and
$K = \I(G\backslash e)$ give a splitting of $\I(G)$.
\end{definition}

The following theorem characterizes all splitting edges in a simple 
graph.  Recall
that $N(u)$ denotes the set of distinct neighbors of the vertex $u$.

\begin{theorem}\cite[Theorem 3.4]{HVT} \label{classify splitting 
edges}
An edge $e = uv$ is a splitting edge of $G$
if and only if $N(u) \subseteq (N(v) \cup \{v\})$ or 
$N(v) \subseteq (N(u) \cup \{u\})$.  
\end{theorem} 

\begin{proof} We shall sketch the main ideas of the proof. 
Let $N(u) \backslash \{v\} = \{u_1, \dots, u_n\}$ and $N(v) \backslash 
\{u\} = \{v_1, \dots, v_m\}$. 
Let $H = G \backslash (N(u) \cup N(v))$. It can be seen that
\begin{align*}
\G(J \cap K) = & \{uvu_i ~|~ u_i \not\in (N(u) \cap N(v))\} \cup 
\{uvv_j ~|~ v_j \not\in (N(u) \cap N(v))\} \ \cup \\
& \{uvz ~|~ z \in (N(u) \cap N(v))\} \cup \{uvm ~|~ m \in \I(H)\}.
\end{align*}

To prove the ``if'' direction, we observe that if $N(u) \subseteq (N(v) 
\cup \{v\})$, then
$$\G(J \cap K) = \{uvv_j ~|~ j = 1, \dots, m\} \cup \{uvm ~|~ m \in 
\I(H)\}.$$
This allows us to construct a splitting function $\G(J \cap K) 
\rightarrow \G(J) \times \G(K)$ as follows
$$w \mapsto (\phi(w), \varphi(w)) = \left\{ \begin{array}{lll} 
(uv,vv_j) & \mbox{if} & w = uvv_j \\ 
(uv,m) & \mbox{if} & w = uvm. \end{array} \right.$$

The ``only if'' direction is proved by proving the contrapositive. Assume that
 $N(u) \not\subseteq (N(v) \cup \{v\})$ and $N(v) \not\subseteq (N(u) 
\cup \{u\})$.
 Then there exist vertices $x,y \in V_G$ such that $ux,vy \in E_G$ and 
$uy,vx \not\in E_G$. 
Now, by using Definition \ref{defn: split} we can show that there 
does not exist a splitting function
 $\G(J \cap K) \rightarrow \G(J) \times \G(K)$.
\end{proof}

\begin{example}\label{splitting edge example} (Splitting edges) 
Consider the following graph $G$:
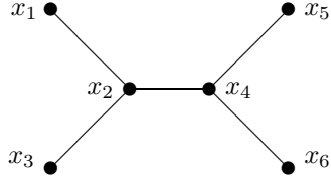
\[
\begin{picture}(100,80)(0,-35)
\put(0,0){\circle*{5}}
\put(-16,-2){$x_2$}
\put(30,0){\circle*{5}}
\put(36,-2){$x_4$}
\put(60,30){\circle*{5}}
\put(66,28){$x_5$}
\put(60,-30){\circle*{5}}
\put(66,-28){$x_6$}
\put(-30,-30){\circle*{5}}
\put(-46,-28){$x_3$}
\put(0,0){\line(1,0){30}}
\put(-30,30){\circle*{5}}
\put(-45,28){$x_1$}
\put(-30,30){\line(1,-1){30}}
\put(-30,-30){\line(1,1){30}}
\put(30,0){\line(1,1){30}}
\put(30,0){\line(1,-1){30}}
\end{picture}\]
The edge $x_1x_2$ is a splitting edge of $G$, but $x_2x_4$ is not.
\end{example}

Once one has identified splitting edges, one can apply the following 
formula.

\begin{theorem}\cite[Theorem 3.6]{HVT} \label{recursive formula 3}
Let $e = uv$ be a splitting edge of $G$, and set $H = G \backslash 
(N(u) \cup N(v))$.  
Set $n = |N(u) \cup N(v)| - 2$.
Then for all $i \geq 1$ and all $j \geq 0$  
\[\beta_{i,j}(\I(G)) =  \beta_{i,j}(\I(G \backslash e)) + 
\sum_{l=0}^i \binom{n}{l} \beta_{i-1-l,j-2-l}(\I(H))\]
where $\beta_{-1,0}(\I(H)) = 1$ and $\beta_{-1,j}(\I(H))=0$ if $j > 
0$.
\end{theorem}

\begin{proof} By Theorem \ref{classify splitting edges} we must have 
$N(u) \subseteq (N(v) \cup \{v\})$ or $N(v) \subseteq (N(u) \cup 
\{u\})$. Suppose 
$N(u) \subseteq (N(v) \cup \{v\})$ and let $N(v) \backslash \{u\} = 
\{v_1, \dots, v_n\}$. 
As before, observe that
$$J \cap K = uv((v_1, \dots, v_n) + \I(H)).$$
The statement now follows by combining Theorem \ref{prop: ekf} and the 
observation that the 
resolution of $(v_1, \dots, v_n) + \I(H)$ can be derived from the 
tensor product of the resolutions of 
$R/(v_1, \dots, v_n)$ and $R/\I(H)$ (since $H$ does not contain any 
vertices 
from among $\{v_1, \dots, v_n\}$). 
\end{proof}

Theorem \ref{recursive formula 3} allows us 
to recover most of the known results in the literature about
the graded Betti numbers of an edge ideal of a forest with fuller generality.  We begin 
by using our
formula to give a recursive formula for the graded Betti numbers of $\I(G)$ when
$G$ is a forest.  This
formula was first proved in \cite{J,JK} via different means.
In fact, our result is slightly more general since it applies to any 
leaf
of a forest, while  \cite{J,JK} required that  a special leaf  be 
removed.  Recall
that by a leaf of $G$, we are referring to an edge with a vertex of 
degree 1.  

\begin{corollary}\label{treeformula}
Let $e=uv$ be any leaf of a forest $G$.
If $\deg v =n$ and $N(v) = \{u,v_1,\ldots,v_{n-1}\}$, then for $i \geq 
1$ and $j \geq 0$ 
\[\beta_{i,j}(\I(G)) = \beta_{i,j}(\I(T)) 
+ \sum_{l=0}^i \binom{n-1}{l} \beta_{i-1-l,j-2-l}(\I(H))\]
where $T = G \backslash e = G \backslash \{u\}$ and $H = G \backslash
\{u,v,v_1,\ldots,v_{n-1}\}$.
Here $\beta_{-1,0}(\I(H)) = 1$ and $\beta_{-1,j}(\I(H))=0$ if $j > 0$.
\end{corollary}

\begin{proof}
The hypotheses imply that $\deg u =1$.  Since $N(u) \subseteq (N(v) 
\cup \{v\})$,
$uv$ is a splitting edge.  Now   
apply Theorem \ref{recursive formula 3}.  The formula is recursive 
since $T$ and $H$
are forests and so they each have a leaf, that is, a splitting edge.
\end{proof}

\begin{remark} In \cite{HVT2} the authors recently showed
that there is in fact a recursive formula to compute the
graded Betti numbers for all chordal graphs.  To prove this
result,
one must schow that a chordal graph always has at least one splitting
edge.
\end{remark}

Theorem \ref{recursive formula 3} can also be used to relate algebraic 
invariants, such as the 
regularity and the projective dimension, of an edge ideal of a graph to 
that of subgraphs.

\begin{corollary}\cite[Corollary 3.7]{HVT} \label{removing edge}
Let $e = uv$ be a splitting edge of a graph $G$, and let $H = G 
\backslash (N(u) \cup N(v))$. Let 
$n = |N(u) \cup N(v)| - 2$. Then we have
\begin{enumerate}
\item $\reg(\I(G)) = \max \{ 2, \reg(\I(G \backslash e)), \reg(\I(H))+1 
\}.$
\item $\pdim(\I(G)) = \max \{ \pdim(\I(G \backslash e)), 
\pdim(\I(H))+n+1 \}.$
\end{enumerate}
\end{corollary}
\noindent Notice that Corollary \ref{removing edge}(2) generalizes 
\cite[Theorem 4.8]{JK} 
which proved the formula in the case that $G$ was a forest.

\begin{remark}
Besides recovering the main theorem of \cite{JK}, Theorem 
\ref{recursive formula 3} 
enables us to recover Theorem \ref{tree reg} 
in the case of forests as first proved by Zheng 
\cite{Z}.
Precisely, we can produce a new 
combinatorial proof of Theorem \ref{tree reg} by
using induction on the number of edges in the graph and combining 
Theorem \ref{recursive formula 3} and Corollary \ref{removing edge} 
(see \cite[Corollary 3.11]{HVT} for
the details). 
\end{remark}

\subsection{Splitting vertices}

Let $G$ be a simple graph, and let $v$ be a vertex of $G$
with $N(v) = \{v_1,\ldots,v_d\}$.  This
section complements the results of the previous section by determining
when $\I(G) = J + K$ with $J =  (vv_1,\ldots,vv_d)$ and $K = 
\I(G\backslash\{v\})$
is a splitting of $\I(G)$.

\begin{lemma} \label{betti intersection 2}
With the notation as above, set
\begin{eqnarray*}
G_i & := & G \backslash (N(v) \cup N(v_i)) ~\mbox{for $i = 1,\ldots,d,$ 
and } \\ 
G_{(v)} & := & G_{\{v_1,\ldots,v_d\}} \cup \{e\in E_G ~|~ \mbox{$e$ 
incident to $v_1,\ldots,v_d$, 
but not $v$}\}.
\end{eqnarray*}
Then
\[J \cap K = v\I(G_{(v)}) + vv_1\I(G_1) + vv_2\I(G_2)+\cdots + 
vv_d\I(G_d).\]
\end{lemma}

\begin{example} Consider the same graph $G$ as in Example 
\ref{splitting edge example}. Take $v = x_2$. Then 
$N(v) = \{v_1=x_1, v_2 = x_3, v_3 = x_4\}$. 
In this case, $G_1 = G_2$ is the graph with two isolated vertices $x_5$ 
and $x_6$, $G_3$
is the empty graph, and $G_{(v)}$ is the graph 
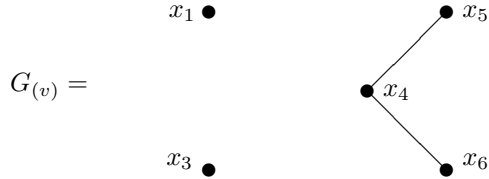
\[ \begin{picture}(200,80)(0,-35)
\put(-5,0){$G_{(v)} = $}
\put(130,0){\circle*{5}}
\put(136,-2){$x_4$}
\put(160,30){\circle*{5}}
\put(166,28){$x_5$}
\put(160,-30){\circle*{5}}
\put(166,-28){$x_6$}
\put(70,-30){\circle*{5}}
\put(54,-28){$x_3$}
\put(70,30){\circle*{5}}
\put(55,28){$x_1$}
%\put(-30,30){\line(1,-1){30}}
%\put(-30,-30){\line(1,1){30}}
\put(130,0){\line(1,1){30}}
\put(130,0){\line(1,-1){30}}
\end{picture} \]
\end{example}

Observe that if $v \in V_G$ is such that $\deg v = 0$, then
$\beta_{i,j}(\I(G)) 
= \beta_{i,j}(\I(G \backslash \{v\}))$ for all $i,j\geq 0$
since $\I(G) = \I(G\backslash\{v\})$.  It therefore
suffices to compute the graded Betti numbers of the edge ideals of each
connected component
of $G$ that contains one or more edges.  If $\deg v = d > 0$ 
and if $G \backslash\{v\}$ consists of isolated vertices,  
then $G = \K_{1,d}$, the complete bipartite graph of
size $1,d$. In this situation, the graded Betti numbers of $\I(G)$ are 
completely known
as seen in Example \ref{k1dbettinumbers}. 
We now give a name to a vertex $v \in V_G$ that does not fall into
either of the above two cases.

\begin{definition} A vertex $v \in V_G$ is a {\it splitting vertex}
if $\deg v = d > 0$ and $G \backslash \{v\}$ is not the graph
of isolated vertices.
\end{definition}

Our choice of name is suitable 
in light of the following theorem.

\begin{theorem}\cite[Theorem 4.2]{HVT}\label{prop: edgeidealssplit}
Let $v$ be a splitting vertex of $G$ with $N(v) = \{v_1,\ldots,v_d\}$.  
Then $\I(G)$ is 
a splittable monomial ideal with a splitting given by $J = 
(vv_1,\ldots,vv_d)$ and 
$K = \I(G \backslash \{v\})$.
\end{theorem}

\begin{proof} By Lemma \ref{betti intersection 2} we have an explicit 
description for 
the generators of $\G(J \cap K)$. A splitting function 
$\G(J \cap K) \rightarrow \G(J) \times \G(K)$ is then given by $w 
\mapsto (\phi(w), \varphi(w))$ where
$$\phi(w) = \left\{ \begin{array}{lll} vv_i & \mbox{if} & w = vv_iv_j 
\in \G(v\I(G_{(v)})) \mbox{ and } i < j \\
vv_i & \mbox{if} & w = vv_iy \in \G(v\I(G_{(v)})) \mbox{ and } y 
\not\in N(v) \\
 vv_i & \mbox{if} & w = vv_im \in \G(vv_i\I(G_i)) \end{array} 
\right.$$
and
$$\varphi(w)= \left\{ \begin{array}{lll} v_iv_j & \mbox{if} & w = 
vv_iv_j \in \G(v\I(G_{(v)})) \mbox{ and } i < j \\
v_iy & \mbox{if} & w = vv_iy \in \G(v\I(G_{(v)})) \mbox{ and } y 
\not\in N(v) \\ 
m & \mbox{if} & w = vv_im \in \G(vv_i\I(G_i)). \end{array} \right.$$
\end{proof}

Applying Theorem \ref{prop: ekf} we obtain:

\begin{theorem}\cite[Theorem 4.6]{HVT} \label{recursive formula 2}
Let $v$ be a splitting vertex of $G$ with $N(v) = \{v_1,\ldots,v_d\}$. 
Let $G_{(v)}$ and 
$G_i$ for $i=1, \dots, d$ be defined as in Lemma \ref{betti 
intersection 2}. Then
\[\beta_{i,j}(\I(G)) = \beta_{i,j}(\I(\K_{1,d})) + \beta_{i,j}(\I(G 
\backslash \{v\})) + 
\beta_{i-1,j}(L)\]
where  $L = v\I(G_{(v)}) + vv_1\I(G_1) + \cdots +vv_d\I(G_d)$ and 
$\K_{1,d}$ is the 
complete bipartite graph of size $1,d$.
\end{theorem}

\begin{proof} One combines Lemma \ref{betti intersection 2}, 
Theorem \ref{prop: edgeidealssplit} and Theorem \ref{prop: ekf}.
\end{proof}

Theorem \ref{recursive formula 2} allows us to relate algebraic 
invariants, such as the 
regularity and the projective dimension, of an edge ideal of a graph to 
that of subgraphs. 

\begin{corollary}\cite[Corollary 4.4]{HVT} \label{removing vertex}
Let $v \in V_G$ be any vertex of a simple graph $G$, and suppose $\deg 
v = d$. Then
\begin{enumerate}
\item $\reg(\I(G)) \ge \max \{ 2, \reg(\I(G \backslash \{v\})) \}.$
\item $\pdim(\I(G)) \ge \max \{d-1, \pdim(\I(G \backslash \{v\}))\}.$
\end{enumerate}
\end{corollary}
\noindent 
Jacques \cite[Proposition 2.1.4]{J} first proved
Corollary \ref{removing vertex} (2)  in the case that $v$ is a {\it 
terminal vertex}, i.e. when $\deg v \le 1$.

Unlike Theorem \ref{recursive formula 3}, we cannot extract a recursive 
formula from
Theorem \ref{recursive formula 2} since the ideal $L =  v\I(G_{(v)}) + 
vv_1\I(G_1) + \cdots +vv_d\I(G_d)$
is not the edge ideal of a graph.  However, a recursive formula
for the graded Betti numbers in the linear strand can be deduced from
Theorem \ref{recursive formula 2}:
\begin{corollary}  Let $v$ be a splitting vertex of a graph $G$.  Then
for all $i \geq 0$,
\[\beta_{i,i+2}(\I(G)) = \beta_{i,i+2}(\I(\K_{1,d})) + 
\beta_{i,i+2}(\I(G\backslash \{v\}))
+\beta_{i-1,i+1}(\I(G_{(v)})).\]
\end{corollary}
\begin{proof}
By evaluating the formula of Theorem \ref{recursive formula 2} at $j = 
i+2$ we get
\[\beta_{i,i+2}(\I(G)) = \beta_{i,i+2}(\I(\K_{1,d})) + 
\beta_{i,i+2}(\I(G\backslash \{v\}))
+\beta_{i-1,i+2}(L).\]
Since $vv_1\I(G_1) + \cdots + vv_d\I(G_d)$ is generated by monomials of 
degree 4 and
$v\I(G_{(v)})$ is generated by monomials of degree 3, we have 
$\beta_{i-1,i+2}(L) = 
\beta_{i-1,i+2}(v\I(G_{(v)})) = \beta_{i-1,i+1}(\I(G_{(v)}))$.
\end{proof}

Theorem \ref{recursive formula 2} also enables us to give new 
combinatorial proofs for
many interesting results.  For example,  the characterization of 
property $N_{2,p}$
of Eisenbud et. al.
for quadratic square-free monomial ideals can be proved without the use 
of
Hochster's formula.  We have included a sketch of
the new proof (see \cite[Corollary 4.7]{HVT} for complete details).

\begin{corollary}\cite[Theorem 2.1]{eghp} \label{N2P2}
Let $G$ be a simple graph with edge ideal $\I(G)$. Then $\I(G)$ 
satisfies property $N_{2,p}$ with $p > 1$ 
if and only if every minimal cycle in $G^c$ has length $\ge p+3$.
\end{corollary}

\begin{proof} One can prove this statement by induction on $|V_G|$, the 
number 
of vertices.  By Theorem \ref{recursive formula 2} we have
\[\beta_{i,j}(\I(G)) = \beta_{i,j}(\I(\K_{1,d})) + 
\beta_{i,j}(\I(G\backslash \{v\}))
+\beta_{i-1,j}(L).\]
So, $\I(G)$ satisfies property $N_{2,p}$ if and only if $\I(\K_{1,d})$ 
and $\I(G\backslash \{v\})$
satisfy property $N_{2,p}$
and $L$  satisfies property $N_{3,p-1}$.  By Example 
\ref{k1dbettinumbers}, 
the ideal $\I(\K_{1,d})$ always has a linear resolution, so it 
satisfies property $N_{2,p}$.
The induction hypothesis allows us to show $\I(G\backslash \{v\})$ 
satisfies property $N_{2,p}$.
Finally, $L$ satisfies property $N_{3,p-1}$ if and only if 
$L=v\I(G_{(v)})$ and
$\I(G_{(v)})$ satisfies property property $N_{2,p-1}$.  The heart of
the proof is to then verify that $L = v\I(G_{(v)})$ and $\I(G_{(v)})$ 
satisfies
property $N_{2,p-1}$ if and only if every minimal cycle of $G^c$ 
containing $v$
has length $\geq p+3$.
\end{proof}

\begin{remark}
Fr\"oberg's 
result about edge ideals with linear resolutions (Theorem \ref{linear 
resolution}), and 
the formula of Theorem \ref{prop: linearstrand_no_c4} first proved by 
Roth and 
the second author  are also
examples of corollaries of Theorem \ref{recursive formula 2}.  See 
\cite{HVT} for
the details of these proofs.
\end{remark}

\subsection{Splitting facets}  In this section we will show that the 
notion of 
splittable ideals can also be used quite profitably to study 
resolutions of facet ideals 
in general. Throughout this section $\Delta$ will denote a simplicial 
complex on a vertex set $V_\Delta$.

\begin{definition} \label{3conn} 
Let $\Delta$ be a simplicial complex, and let $F$ be a facet of 
$\Delta$. 
The {\it connected component} of $F$ in $\Delta$, denoted by 
$\lk_\Delta(F)$, 
is the connected component of $\Delta$ containing $F$. 
If $\lk_\Delta(F) \backslash F = \langle G_1, \dots, G_p \rangle$, then 
we define 
the {\it reduced connected component} of $F$ in $\Delta$, denoted by 
$\rlk_{\Delta}(F)$, to 
be the simplicial complex whose facets are given by 
$G_1 \backslash F, \dots, G_p \backslash F$, where if there exist $G_i$ 
and $G_j$ such 
that $\emptyset \not= G_i \backslash F \subseteq G_j \backslash F$, 
then we shall disregard 
the bigger facet $G_j \backslash F$ in $\rlk_\Delta(F)$.
\end{definition}

\begin{example}
Consider the simplicial complex $\Delta$ with the facet set 
$\F(\Delta) = \big\{ \{1,2,3\}, \{1,3,4\}, \{1,4,5\}, \{1,5,6\} 
\big\}$. Let $F = \{1,5,6\}$. 
Then $\lk_\Delta(F) = \Delta$ and $\rlk_\Delta(F)$ is the simplicial 
complex with the facet set 
$\big\{ \{2,3\}, \{4\} \big\}$. Note that $\{1,3,4\} \backslash F = 
\{3,4\}$ contains 
$\{4\} = \{1,4,5\} \backslash F$ so we disregard the bigger set (which 
is $\{3,4\}$) in obtaining $\rlk_\Delta(F)$.
\end{example}

For a facet $F$ of a simplicial complex $\Delta$, we denote by $\Delta' 
= \Delta \backslash F$ the simplicial 
complex obtained by removing $F$ from the facet set of $\Delta$. Let
$$ J = (x^F) \ \text{and} \ K = \I(\Delta').$$
Note that $\G(\I(\Delta))$ is the disjoint union of $\G(J)$ and 
$\G(K)$. We are interested in 
finding $F$ such that $\I(\Delta) = J + K$ gives a splitting of
$\I(\Delta)$.

\begin{definition} \label{defn: splitting facet}
With the above notation, we call $F$ a {\it splitting facet} of 
$\Delta$ if $\I(\Delta) = J + K$ 
is a splitting of $\I(\Delta)$.
\end{definition}

The following result gives a recursive-type formula for the graded 
Betti numbers 
of the facet ideal of a simplicial complex in terms of the Betti 
numbers 
of facet ideals of subcomplexes.  This
result also generalizes Theorem \ref{recursive formula 3}.

\begin{theorem}\cite[Theorem 5.5]{HVT} \label{simplicial:betti}
Let $F$ be a splitting facet of a simplicial complex $\Delta$. Then for 
all $i \geq 1$
and $j \geq 0$ 
$$\beta_{i,j}(\I(\Delta)) = \beta_{i,j}(\I(\Delta')) + \sum_{l_1=0}^i 
\sum_{l_2=0}^{j-|F|} \beta_{l_1-1,l_2}(\I(\rlk_\Delta(F))) 
\beta_{i-l_1-1,j-|F|-l_2}(\I(\Omega))$$
where $\Delta' = \Delta \backslash F$ and $\Omega = \Delta \backslash 
\lk_\Delta(F).$
Here  $\beta_{-1,0}(I) = 1$ and $\beta_{-1,j}(I) = 0$ if $j > 0$ for 
$I = \I(\rlk_\Delta(F))$ and $\I(\Omega)$.
\end{theorem}

\begin{proof} It can be shown that $J \cap K = x^F(\I(\rlk_\Delta(F)) + 
\I(\Omega))$. 
Observe that $\Omega$ and $\rlk_\Delta(F)$, by definition, do not share 
any common vertices. 
Thus, the minimal free resolution of $\I(\rlk_\Delta(F)) + \I(\Omega)$ 
can be derived from the 
tensor product of the resolutions of $\I(\rlk_\Delta(F))$ and 
$\I(\Omega)$. The result now follows
 by applying Theorem \ref{prop: ekf}.
\end{proof}

We will now show that our formula in Theorem \ref{simplicial:betti} is 
recursive when $\Delta$ is a 
simplicial forest.
To do so, we first show that a leaf of $\Delta$ is a splitting facet.
Recall that if $F$ is a leaf of $\Delta$, then $F$ must have a vertex 
that does not belong to
any other facet of the simplicial complex (see \cite[Remark 
2.3]{faridi:2002}). 

\begin{theorem}\cite[Theorem 5.6]{HVT} \label{3split}
If $F$ is a leaf of $\Delta$, then $F$ is a splitting facet of 
$\Delta$.
\end{theorem}

\begin{proof} The proof uses a similar line of reasoning as that of the 
``if'' direction of 
Theorem \ref{classify splitting edges}. An explicit description for the 
generators of $\G(J \cap K)$ 
can be given and a splitting function $s: \G(J \cap K) \rightarrow 
\G(J) \times \G(K)$ is constructed
 in the most natural way. The fact that $F$ contains a vertex $x$ that 
does not belong to any other 
facets of $\Delta$ guarantees that the function $s$ does indeed satisfy 
all the conditions
of Definition \ref{defn: split}.
\end{proof}

Recall that a forest is a simplicial complex with the property that 
every 
nonempty connected subcomplex has a 
leaf. Since $\lk_\Delta(F)$, $\Delta \backslash F$ and $\Delta 
\backslash \lk_\Delta(F)$ are 
subcomplexes of $\Delta$, it follows directly from the definition that 
if $\Delta$ is a forest 
then so are $\lk_\Delta(F)$, $\Delta \backslash F$ and $\Delta 
\backslash \lk_\Delta(F)$. 
Thus, to show that our formula in Theorem \ref{simplicial:betti} is 
recursive when $\Delta$ is a 
simplicial forest,
we need to show that $\rlk_\Delta(F)$ is also a forest; this is
the content of the next lemma.

\begin{lemma}\cite[Lemma 5.7]{HVT} \label{3forest}
Let $F$ be a facet of a forest $\Delta$. Then $\rlk_\Delta(F)$ is a 
forest.
\end{lemma}

Our recursive formula for simplicial trees generalizes
Corollary \ref{treeformula} to higher dimensions.

\begin{theorem}\cite[Theorem 5.8]{HVT} \label{simplicial:recursive}
Let $F$ be a leaf of a simplicial forest $\Delta$, and let 
$\Delta' = \Delta \backslash F$ and $\Omega = \Delta \backslash 
\lk_\Delta(F).$ Then
$\Delta', \Omega$, and $\rlk_\Delta(F)$ are also simplicial forests and 
for 
all $i \geq 1$ and $j \geq 0$ 
$$\beta_{i,j}(\I(\Delta)) = \beta_{i,j}(\I(\Delta')) + \sum_{l_1=0}^i 
\sum_{l_2=0}^{j-|F|} \beta_{l_1-1,l_2}(\I(\rlk_\Delta(F))) 
\beta_{i-l_1-1,j-|F|-l_2}(\I(\Omega))$$
where  $\beta_{-1,0}(I) = 1$ and $\beta_{-1,j}(I) = 0$ for $j > 0$ if  
$I = \I(\rlk_\Delta(F))$ or $\I(\Omega)$.
\end{theorem}
Recall that a simplicial complex $\Delta$ is said to be a pure 
$(d-1)$-dimensional simplicial complex
if 
$\dim F = d-1$, i.e., $|F| = d$, for any facet $F$ of $\Delta$. For a 
face $G$ of 
dimension $d-2$ of a pure $(d-1)$-dimensional simplicial complex 
$\Delta$ we define the 
{\it degree} of $G$, written $\deg_\Delta(G)$, to be the cardinality of 
the set
 $\{ F \in \F(\Delta) \ | \ G \subseteq F \}$. Let $\A(\Delta)$ denote 
the set of
 $(d-2)$-dimensional faces of $\Delta$. The following result gives a 
formula for the graded Betti numbers 
in the linear strand of the facet ideal of a pure simplicial forest. 

\begin{theorem}\cite[Theorem 5.9]{HVT} \label{3linear-strand}
Let $\Delta$ be a pure $(d-1)$-dimensional forest (for some $d \ge 2$). 
Then
$$\beta_{i,i+d}(\I(\Delta)) = \left\{ \begin{array}{lcl} |\F(\Delta)| & 
\text{if} & i = 0 \\
{\displaystyle \sum_{G \in \A(\Delta)} {\deg_\Delta(G) \choose i+1}} & 
\text{if} & i \ge 1. 
\end{array} \right.$$
\end{theorem}

\begin{proof} The proof uses induction on the number of facets of 
$\Delta$ and 
makes use of Theorem \ref{simplicial:recursive}.
\end{proof}
\noindent Theorem \ref{3linear-strand} was first proved by Zheng 
\cite{Z} under 
the extra condition that $\Delta$ is connected in codimension 1. By 
using the notion of splittable 
ideals, this hypothesis can be removed.  When $d=2$, then $\Delta$
is a forest in the standard sense, and we recover Corollary 
\ref{foreststrand}.  We
can therefore view Theorem \ref{3linear-strand} as a partial 
generalization
of Theorem \ref{prop: linearstrand_no_c4}.

%%%%%%%%%%%%%%%%%%%%%%%%%%%%%%%%%%%%%%%%%%%%%%%%%%%%%%

\section{Open questions}

As noted in the introduction, one of our goals in writing this survey 
is to promote further
research on the resolutions of square-free monomial ideals from a facet 
ideal
point-of-view.  We end this paper with
some natural questions whose answers we would be interested in 
knowing.

\subsection{Building a dictionary} One of the themes stressed in this 
paper
is how the combinatorial data of either a graph $G$ 
or a simplicial complex $\Delta$ appears in the minimal graded free 
resolution
of $\I(G)$ or $\I(\Delta)$.  Although we have shown that many graded 
Betti
numbers can be described directly from the combinatorial data, there is 
still much
we do not know.  We therefore pose the general question:

\begin{question}
Let $\I(\Delta)$ be the facet ideal of a simplicial complex. For which 
$i$ and $j$ is 
there a formula for the graded Betti number $\beta_{i,j}(\I(\Delta))$ in 
terms of combinatorial data of $\Delta$?
\end{question}

This question is probably too imprecise; one should not expect a simple 
answer to
this question because some of the numbers $\beta_{i,j}(\I(\Delta))$ 
will
depend upon the characteristic of the field (as seen in Reisner's 
example
in the introduction).
However, an interesting place to start is to see what structure 
$\Delta$ must have to force 
$\beta_{i,j}(\I(\Delta))$ to be zero.

Among the graded Betti numbers of a facet ideal, those 
in the linear strand are of particular interest. We raise the following 
question:

\begin{question} \label{linear strand question}
Let $\Delta$ be a pure $(d-1)$-dimensional simplicial complex. Is there 
a formula for $\beta_{i,i+d}(\I(\Delta))$ which describes the linear 
strand of the resolution of $\I(\Delta))$ similar to that of Theorem 
\ref{3linear-strand} and Corollary \ref{prop: linearstrand_no_c4}?
\end{question}

When the simplicial complex is of dimension one (i.e., it is a graph) 
Question \ref{linear strand question} has been addressed positively by 
Roth and the second author (see Corollary \ref{prop: linearstrand_no_c4}) 
for edge ideals of graphs having no minimal 4-cycles. It is natural to 
ask the question for graphs which contain minimal 4-cycles.

\begin{question}
If $G$ has minimal 4-cycles, what is a formula for 
$\beta_{i,i+2}(\I(G))$?
\end{question}

Finally, when $\I(G)$ is an edge ideal, we know $\beta_{0,j}(\I(G))$ 
and $\beta_{1,j}(\I(G))$
and $\beta_{2,j}(\I(G))$ can be computed 
directly from the graph (see Examples \ref{beta0j}, \ref{beta1j}, and 
\ref{beta2j}, respectively) for all $j$.  
Hence, the first place to search for new formulas for 
$\beta_{i,j}(\I(G))$ is in the case that $i = 3$.
Note that one could take the approach used in Example \ref{beta25} by 
compiling a list of all graphs
$H$ with $\dim_k \widetilde{H}_i(\Delta(H),k) > 0$.  More compact 
formulas, however, would be preferred.  Indeed,
it would be nice to find an alternative formula for 
$\beta_{2,5}(\I(G))$ that avoided having to identifying large
numbers of induced subgraphs.

Of course, once we have built a reasonably sized dictionary, we want to 
use the
tools of commutative algebra to answer questions about graph theory.  
For example,
Roth and the second author \cite{RVT} showed that if one uses the 
Bigatti-Hullet-Pardue
theorem about the growth of graded Betti numbers of lexicographical 
ideals, one
can obtain a crude bound on the number of triangles in a graph.  Do 
other
such results await us?

\subsection{Characteristic-independence}  
These questions are inspired by Katzman's
paper \cite{K}.  From Hochster's formula 
(see Theorem \ref{prop:  hochster}) it follows that the 
graded Betti numbers of a monomial ideal may depend upon the field $k$.  
Reisner's example
(see Section 4 of \cite{K}) 
of the triangulation of the real projective plane is a classical 
example of how the graded
Betti numbers depend upon  char($k$).
Highlighted below are some of Katzman's results on  how the
numbers $\beta_{i,j}(\I(G))$ of an edge ideal $\I(G)$ depend upon 
char$(k)$.
\begin{theorem}\label{char depend}
Let $\I(G)$ be the edge ideal of a simple graph $G$.  Then
\begin{enumerate}
\item[$(i)$] $\beta_{i,j}(\I(G))$ is independent of char$(k)$ for all 
$i \leq 5$.
\item[$(ii)$] $\beta_{i,j}(\I(G))$ is independent of char$(k)$ if 
$|V_G| \leq 10$.
\item[$(iii)$] there exists exactly 4 non-isomorphic graphs $G$ with 
$|V_G| = 11$ such that 
the numbers $\beta_{i,j}(\I(G))$ depend upon char$(k)$.  In each case
the Betti number $\beta_{i,j}(\I(G))$ that depends upon char$(k)$ has 
$i = 7$ or $8$,
and char$(k) = 2$.
\end{enumerate}
\end{theorem}
Theorem \ref{char depend} (i) extends an earlier result of Terai and 
Hibi \cite{TH}
that $\beta_{2,j}(\I(G))$ and $\beta_{3,j}(\I(G))$ are independent of 
char$(k)$.
The above theorem does not tell us whether the numbers 
$\beta_{6,j}(\I(G))$
depend upon char$(k)$. 
By Theorem \ref{thm: betti edge range} we know that $\beta_{6,j}(\I(G)) 
\neq 0$
only if $j = 8,\ldots 14$.  Katzman was able to show that 
$\beta_{6,j}(\I(G))$
does not depend upon the characteristic if $j \neq 12$, but left
open the question when $j=12$.  This brings us to our first question.

\begin{question}
Does the number $\beta_{6,12}(\I(G))$ depend upon 
char$(k)$?
\end{question}

Theorem \ref{char depend} (iii) says that one must consider graphs $G$ 
with at least
12 vertices to find an example.  
Related to this question, we could ask:

\begin{question}
Can we identify graphs $G$ (or simplicial complexes $\Delta$) through 
some combinatorial means (e.g., $G$ or $\Delta$ has a subgraph or a 
subcomplex of a particular form) with
the property that the graded Betti numbers of their edge ideals depend 
on char$(k)$?
\end{question}

An answer to this more general question might then provide an easy 
answer to our 
first question.
Note that in the paper of Katzman \cite{K} all of the examples have the 
property that
the numbers $\beta_{i,j}(\I(G))$ only change if char$(k) =2$.  One is 
naturally
lead to ask if the graded Betti numbers only change if  char$(k) = 2$. 
As explained to us by Katzman, the answer to this question is no since 
for any prime $p$, one can always find a graph 
$G$ so that the graded Betti numbers of $\I(G)$ are different in 
char$(k)=p$ and char$(k)=0$.
To find such a graph, begin with any simplicial complex whose homology 
depends upon
the characteristic that you desire.  Then construct the barycentric 
subdivision of
the simplicial complex where the new vertices are the old nonempty 
faces, and the 
new faces are chains of old nonempty faces.  Construct a graph $H$ 
whose
vertices are the old faces and whose edges are pairs of incomparable 
faces.  Then
the homology of the simplicial complex is the same as its barycentric 
subdivision.
Furthermore, the barycentric subdivision is the clique complex 
associated to $H$.  
Because $\I_{\Delta(H)} = \I(H^c)$, the graded Betti numbers of 
$\I(H^c)$
change if the characteristic is $p$.

\subsection{Other algebraic invariants and properties} Besides the 
graded Betti numbers, we are also 
interested in other algebraic invariants and properties of edge ideals 
and facet ideals which are related to
Betti numbers. For example, the regularity, which measures the width of 
the resolution; the projective dimension, 
which measure the length of the resolution; and property $N_{d,p}$, 
which captures how long the resolution must 
have linear syzygies. 

Inspired by Corollary \ref{tree reg}, we raise the following question.

\begin{question} \label{complex reg}
Let $\Delta$ be a simplicial complex. Is there a formula that relates 
$\reg(\I(\Delta))$ to combinatorial data,
e.g., the number of subcomplexes of a particular type, of $\Delta$?
\end{question}

Note that Question \ref{complex reg} has been partially
answered by the two authors \cite{HVT2} in the case of edge
ideals, i.e.,  Theorems \ref{tree reg} and \ref{upperbound} give lower
and upper bounds.
Another direction is to generalize Theorem 
\ref{tree reg} to higher dimension, i.e., to find a formula for 
$\reg(\I(\Delta))$ when $\Delta$ is a pure simplicial tree.   Some
partial results in this direction can also be found in \cite{HVT2}.

Since the method of mathematical induction has proved to be significant 
in obtaining many of our results, it would be interesting to relate 
algebraic invariants and properties of edge ideals (or facet ideals) of 
graphs (or simplicial complexes) to that of subgraphs (or subcomplexes). 
We propose to seek for generalizations of Corollaries \ref{removing 
edge} and \ref{removing vertex}.

\begin{question}
Let $\Delta$ be a simplicial complex. Let $F$ be a facet of $\Delta$ 
and let $v$ be a vertex of $\Delta$.
\begin{enumerate}
\item Is there a formula which relates $\reg(\I(\Delta))$ and 
$\pdim(\I(\Delta))$ to $\reg(\I(\Delta \backslash F))$ and $\pdim(\I(\Delta 
\backslash F))$?
\item Is there a formula which relates $\reg(\I(\Delta))$ and 
$\pdim(\I(\Delta))$ to $\reg(\I(\Delta \backslash \{v\}))$ and $\pdim(\I(\Delta 
\backslash \{v\}))$?
\end{enumerate}
Here, by $\Delta \backslash \{v\}$ we mean the simplicial complex one 
obtains by removing
from  $\Delta$ all facets that contain $v$.
\end{question}

In connection to Green's famous conjecture on canonical curves, the 
property that an ideal has 
linear syzygies up to a given step has sparked much research. One of 
these is the characterization 
for property $N_{2,p}$ of edge ideals due to Eisenbud, Green, Hulek and 
Popescu (see Corollary \ref{np property}). 
We would like to find a similar characterization in higher dimension.

\begin{question}
Let $\Delta$ be a pure $(d-1)$-dimensional simplicial complex. Is there 
a necessary and sufficient condition, based 
upon combinatorial data of $\Delta$, for $\I(\Delta)$ to satisfy 
property $N_{d,p}$ for $p > 1$, i.e., 
$\I(\Delta)$ is generated in degree $d$ (which is obvious) and has 
linear syzygies up to the $p$-th step?
\end{question}

\subsection{Splittable monomial ideals} In the previous section we 
demonstrated the usefulness of the notion
of splittable ideals to study the graded Betti numbers of edge and 
facet ideals.  Our
remaining questions are interested in extending some of these ideas.

In Section 4 we considered two natural splittings of an edge ideal in 
terms
of two natural graph operations, namely, removing an edge and removing 
a vertex. We
can ask if there are any other ways to split the generators of an edge 
ideal.

\begin{question}
If $\I(G)$ is the edge ideal of a graph $G$, are there other splittings 
$\I(G) = J + K$
that give us information on the graded Betti numbers of $\I(G)$ in 
terms of the subgraphs of $G$?
\end{question}

We have classified all splitting edges for edge ideals.  We have also shown 
that for facet ideals,
the facets corresponding to leaves are splitting facets.  However, we 
have
left open the question if there are any other splitting facets.  So one 
is lead to ask:

\begin{question} \label{facet characterization}
Is there a classification of splitting facets?
\end{question}

An affirmative answer to Question \ref{facet characterization} was 
recently obtained 
by the two authors \cite{HVT2}.  This 
classification of splitting facets has lead to interesting 
consequences, 
for example, a formula calculating the regularity of a simplicial 
forest similar to that of Theorem \ref{tree reg}. 

While we introduced a notion of a splitting vertex for graphs, we have 
not
identified an analog of this concept for simplicial complexes. So we 
can ask:

\begin{question} \label{splitting vertex question}
What is the correct generalization of a splitting vertex in the context 
of
simplicial complexes? 
\end{question}

We saw in \cite{HVT} that if $v$ is a non-isolated vertex of a graph 
$G$ such that
 $K = \I(G \backslash \{v\}) \not= (0)$, then $\I(G) = J + K$ is a 
splitting of $\I(G)$ 
(here, $J = ( \{vx ~|~ vx \in E_G \})$). The next example shows that 
the same 
phenomenon is not true for simplicial complexes in general. Thus, the 
question of 
characterizing splitting vertices of simplicial complexes is more 
subtle than that for graphs. 

\begin{example} Consider the simplicial complex $\Delta$ with the facet 
set 
$\big\{ F_1 = \{0,1,2\}, F_2 = \{0,3,4\}, F_3 = \{0,5,6\}, 
G_1 = \{1,2,3\}, G_2 = \{3,4,5\}, G_3 = \{5,6,1\} \big\}$. Consider the 
vertex 0 of 
$\Delta$. Then $F_1, F_2, F_3$ are facets of $\Delta$ containing 0. Let 
$J = (x^{F_1}, x^{F_2}, x^{F_3})$ and $K = (x^{G_1}, x^{G_2}, 
x^{G_3})$. We claim that 
$\I(\Delta) = J + K$ is not a splitting for $\I(\Delta)$. Indeed, 
suppose there exists a 
splitting function $s = (\phi,\varphi): \G(J \cap K) \rightarrow \G(J) 
\times \G(K)$ 
satisfying the two conditions of Definition \ref{defn: split}. Let $L_i 
= F_i \cup G_i$ 
for $i = 1, 2, 3$. Then $S = \{x^{L_1}, x^{L_2}, x^{L_3} \}$ is a 
subset of $\G(J \cap K)$. 
Thus, by definition $\lcm(\phi(S))$ must strictly divide $\lcm(S)$. 
Moreover, it is easy to
 see that $\phi(x^{L_i}) = x^{F_i}$ for $i = 1, 2, 3$. However, we now 
have 
$\lcm(\phi(S)) = x^{F_1 \cup F_2 \cup F_3} = x^{\{0,1,2,3,4,5,6\}} = 
\lcm(S)$, 
a contradiction. This shows that 0 cannot be a splitting vertex of 
$\Delta$.
\end{example}

As we saw in \cite{HVT}, the idea of a splitting vertex is used to 
derive
a recursive formula for the graded Betti numbers in the linear strand 
of edge ideals. This also
enables us to give a combinatorial proof for a result of Eisenbud et. 
al. which characterizes property $N_{2,p}$ for edge ideals. Thus, an 
answer to Question \ref{splitting vertex question} might therefore lead 
one to a characterization of property $N_{d,p}$ for the facet ideal of a 
pure $(d-1)$-dimensional simplicial complex.

\section*{Acknowledgments}
We would like to thank all the organizers involved with the
Midwest Algebra, Geometry and their Interactions Conference (MAGIC) for 
a wonderful conference, and for
their encouragement to write this survey.  
We would also especially like to thank both Jessica Sidman and Moty 
Katzman for reading an earlier draft of this paper and for 
providing invaluable suggestions and improvements.  
The first author is partially supported by the Louisiana Board of Regents Enhancement Grant and the second
author acknowledges the funding received by NSERC while working on this 
project.  We also thank the referee for their comments and suggestions

%%%%%%%%%%%%%%%%%%%%%%%%%%%%%%%%%%%%%%%%%%%%%%%%%%%%%%%

\end{document}